\documentclass[review]{elsarticle}
\usepackage[utf8]{inputenc}

\usepackage{lineno,hyperref}
\modulolinenumbers[5]

\newcommand{\VM}[1]{\mathbf{#1}}
\newcommand{\subs}[1]{\mathcal{#1}}
\newcommand{\norm}[1]{\left \| #1 \right \|}
\newcommand{\Ten}[1]{\mathcal{#1}}
\newcommand{\R}{\mathbb{R}}
\newcommand{\C}{\mathbb{C}}

\newcommand{\xkb}[2][k]{\VM{#2}^{(#1)}}
\DeclareUnicodeCharacter{2212}{-}
\newlength\figW
\usepackage{tikz}

\usepackage{pgfplots}
    \pgfplotsset{compat=1.12}
    \usetikzlibrary{positioning,calc,shapes,arrows}

\pgfplotsset{every tick label/.append style={font=\tiny}}
\newenvironment{customlegend}[1][]{%
    \begingroup
    \csname pgfplots@init@cleared@structures\endcsname
    \pgfplotsset{#1}%
}{%
\csname pgfplots@createlegend\endcsname
\endgroup
}%

\def\addlegendimage{\csname pgfplots@addlegendimage\endcsname}
\pgfplotsset{
    cycle list={%
        {draw=black,mark=star,solid},
        {draw=black, mark=square,solid},%densely dashed},
        {draw=black,mark=+,solid},%dashdotted}, %every mark/.append style={rotate=90},
        {black,mark=o},
        {draw=black, mark=none,solid}}}

\definecolor{mycolor1}{rgb}{0.00000,0.44700,0.74100}%
\definecolor{mycolor2}{rgb}{0.85000,0.32500,0.09800}%
\definecolor{mycolor3}{rgb}{0.92900,0.69400,0.12500}%
\definecolor{mycolor4}{rgb}{0.49400,0.18400,0.55600}%

\usepackage{float}
\usepackage{cool}
\usepackage{subcaption}
\pgfplotsset{compat=1.14}
\usepackage{todonotes}
\usepackage{mathtools}
\usepackage{algorithm,algpseudocode}
\usepackage[margin=0.5in]{geometry}
\usepackage{comment}
\newtheorem{theorem}{Theorem}
\newtheorem{eexample}{Example}
\newtheorem{eproof}{Proof}

\newtheorem{lemma}{Lemma}
\newtheorem{definition}{Definition}

\journal{Journal of Computational and Applied Mathematics}

\begin{document}

\begin{frontmatter}

\title{Tensor-Krylov method for computing eigenvalues of parameter-dependent matrices}
% \tnotetext[mytitlenote]{Fully documented templates are available in the elsarticle package on \href{http://www.ctan.org/tex-archive/macros/latex/contrib/elsarticle}{CTAN}.}

%% Group authors per affiliation:
\author[KUL]{Koen Ruymbeek\corref{cor1}}
\cortext[cor1]{Corresponding author}
\ead{koen.ruymbeek@kuleuven.be}

\author[KUL]{Karl Meerbergen}
\author[KUL]{Wim Michiels}
\address[KUL]{Department of Computer Science, KU Leuven, Celestijnenlaan 200A, 3001 Leuven, Belgium}
% \fntext[myfootnote]{Since 1880.}

%% or include affiliations in footnotes:
% \author[mymainaddress,mysecondaryaddress]{Elsevier Inc}
% \ead[url]{www.elsevier.com}

% \author[mysecondaryaddress]{Global Customer Service\corref{mycorrespondingauthor}}
% \cortext[mycorrespondingauthor]{Corresponding author}
% \ead{support@elsevier.com}

% \address[mymainaddress]{1600 John F Kennedy Boulevard, Philadelphia}
% \address[mysecondaryaddress]{360 Park Avenue South, New York}

\begin{abstract}
In this paper we extend the Residual Arnoldi method for calculating an extreme eigenvalue (e.g. largest real part, dominant,...) to the case where the matrices depend on parameters. The difference between this Arnoldi method and the classical Arnoldi algorithm is that in the former the residual is added to the subspace. We develop a Tensor-Krylov method that applies the Residual Arnoldi method (RA) for a grid of parameter points at the same time. The subspace contains an approximate Krylov space for all these points. Instead of adding the residuals for all parameter values to the subspace we create a low-rank approximation of the matrix consisting of these residuals and add only the column space to the subspace. In order to keep the computations efficient, it is needed to limit the dimension of the subspace and to restart once the subspace has reached the prescribed maximal dimension. The novelty of this approach is twofold. Firstly, we observed that a large error in the low-rank approximations is allowed without slowing down the convergence, which implies that we can do more iterations before restarting. Secondly, we pay particular attention to the way the subspace is restarted, since classical restarting techniques give a too large subspace in our case. We motivate why it is good enough to just keep the approximation of the searched eigenvector. At the end of the paper we extend this algorithm to shift-and-invert Residual Arnoldi method to calculate the eigenvalue close to a shift $\sigma$ for a specific parameter dependency. We provide theoretical results and report numerical experiments. The Matlab code is publicly available.
\end{abstract}


\end{frontmatter}

\linenumbers

\section{Introduction}
% \todo[inline]{Foutje dat moet verbeterd worden na revisie: $f_i$ moet niet separabel zijn !}
In this paper, we consider parameter-dependent eigenvalue problems 
\begin{equation} \label{eqn:VW_A}
\VM A(\VM \omega) \VM x(\VM \omega) = \lambda(\VM \omega) \VM x(\VM \omega), \quad \VM \omega \in \Omega, \quad \VM A: \Omega \rightarrow \R^{n \times n}: \VM \omega \mapsto \sum_{i=1}^m f_i(\VM \omega) \VM A_i
\end{equation}
on a compact parameter set $\Omega \subset \R^d$ and functions $f_i, i=1,2,\hdots,m$ where $m$ is small relative to $n$. We look for an eigenpair $(\lambda, \VM x)$ over $\Omega$ for which $\lambda$ is extreme following some ordering, e.g. we sort by modulus if $\lambda$ is the dominant eigenvalue. We assume that the wanted eigenvalue converges for any $\VM \omega$ in a few iterations of the Arnoldi algorithm. In practice, this means that $\lambda$ should be a well separated extreme eigenvalue of $\VM A(\VM \omega)$. Here we want to develop a numerical method to efficiently calculate $\lambda$ for a large set of points in $\Omega$ as calculating this eigenvalue for each point individually using standard eigenvalue solvers becomes soon infeasible. At the end of the paper, we extend our method to a case for which the eigenvalue closest to a shift $\sigma$ is wanted. 
% where we search for the eigenvalue closest to a shift $\sigma$ when all matrices except one is of rank one.

If $\VM A(\VM \omega)$ is symmetric for all $\VM \omega \in \Omega$, a sampling-based subspace method \cite{ruymbeek2019calculating} and the Successive Constraint method (SCM) \cite{sirkovic2016} were recently proposed. In the latter, the desired eigenvalue is approximated by iteratively computing lower-and upper bounds for each parameter value. Another class of solvers are based on polynomial approximation, see for example \cite{Andreev2012}. This technique is also used for the related problem of calculating statistical moments if the parameter is stochastic, see \cite{Benner2019}. The disadvantage of using polynomials is that they suffer from possible lack of smoothness of the desired eigenvalue function \cite{Fenzi2019}.

The method that we propose is an example of a Tensor-Krylov method. In these methods, the iteration vectors in the Krylov method are represented by a low rank tensor. Tensor-Krylov techniques are already discussed for various problems, for example to solve a large system of non-linear equations \cite{Bouaricha1996, Bouaricha1998} and parametrized linear systems \cite{Kressner2011}. 

There are two important aspects in our method. First the iteration vector is represented by a tensor which is truncated to a low-rank representation, for reasons of storage and computational efficiency. Second, the Ritz pairs are obtained from projection of $\VM A(\VM \omega)$ on the sum of the Krylov spaces for all sample points. A practical issue is the growth of the dimension of this subspace as we add multiple vectors in each iteration. Therefore we often need to restart to keep the algorithm feasible. We gain in efficiency if we need to add less vectors to the subspace (and so introducing a larger error) and have a good way to restart the subspace. We demonstrate that adding fewer vectors in a careful way does not lead to slower convergence and we give an explanation. We give practical bounds for the introduced errors.

%\medskip
%\paragraph{Outline and notation} 
The plan of the paper is as follows. In Section 2 we discuss the Arnoldi method and analyse how errors are propagated in the algorithm and give convergence results. We explain the proposed algorithm in the one parameter case and discuss why it is an extension of the Residual Arnoldi algorithm in Section 3. Subsequently, we illustrate the convergence and error propagation with some numerical experiments. In Section 5, we revisit the method for multiple parameters and report a numerical experiment. We extend our method to shift-and-invert Residual Arnoldi when all matrices except one is of rank one and give some numerical examples. We finish this paper with some concluding remarks. The Matlab code can be downloaded following \href{http://twr.cs.kuleuven.be/research/software/delay-control/TensorKrylov}{this link}.

Throughout the paper, we denote vectors by small bold letters and matrices by capital bold letters. The subspace of $\R^n$ spanned by the columns of an orthonormal matrix $\VM V$ is denoted by the calligraphic letter $\subs{V}$. An index $k$ may be added to refer to the iteration where this subspace is made.
With calligraphic letters other than $\subs{V}$, we denote tensors of order strictly larger than two. Further, we denote by $\VM{I}$ the identity matrix of appropriate dimensions and by the function $\text{diag}(\VM{v})$ a diagonal matrix with the elements of the vector $\VM{v}$ on its diagonal. With $\norm{.}$ we denote the Euclidean norm and we use  $\norm{.}_F$ to refer to the Frobenius norm.

We end this section by introducing the concept of reduced eigenvalue problem. Let $\subs{V} \subset \R^n$ be a subspace and $\VM V$ an orthonormal matrix where the columns form a basis for $\subs{V}$, then we call $\VM{V}^T \VM{A} \VM{V}$ the \emph{reduced eigenvalue problem}. We denote by $\left( \lambda^\subs{V}, \VM{x}^\subs{V} \right)$ an eigenpair of this reduced eigenvalue problem and the corresponding Ritz vector by $\VM y = \VM V \VM x^\subs{V}$. A subindex may be used to indicate the iteration where the subspace is made. We define the corresponding residual as $\VM r = \VM A \VM y - \lambda^\subs{V} \VM y$. We further assume that an eigenvector has norm $1$. All proposed algorithms in this paper are implemented in Matlab version R2018a. All experiments are performed on an Intel i5-6300U with $2.5$ GHz and $8$ GB RAM.

\section{Arnoldi method}
In this section we introduce the Arnoldi algorithm in a general form, we explain how convergence is measured and give some results regarding its performance when the iterations are contaminated by error. The Arnoldi algorithm itself can be found in pseudo-code in Algorithm \ref{algo:Arn}. In the Arnoldi algorithm, a subspace $\subs{V}$ is iteratively built and we look for an approximation of the searched eigenvector in this subspace. In the algorithm, we need to decide about the choice of the continuation vector $\VM t$, see line \ref{line:tk} in Algorithm \ref{algo:Arn}, which, after multiplying by $\VM A \VM V$, is added to the subspace.

\begin{algorithm}[h] 
\hspace*{\algorithmicindent} \textbf{Input:} Matrix $\VM{A}$, normalized vector $\VM{u} \in \R^n$ \\
\hspace*{\algorithmicindent} \textbf{Output:} Approximate eigenpair $(\hat \lambda, \VM{y})$ of $\VM A$ 
\begin{algorithmic}[1]
\State $\VM w^\bot = \VM u$ 
\For{$k=1$ to $n-1$}
\State Add $\VM w^\bot$ to the subspace $\subs V$ 
\State $\VM w = \VM A \VM V \VM t$  \label{line:tk}
\State $\VM w^\bot =  \left( \VM I - \VM V \VM V^T \right) \VM w$
\State Calculate eigenpair  $\left( \lambda^{\subs{V}}, \VM x^{\subs{V}} \right)$ of $\VM V^T \VM A \VM V$
\State Break down if $\left( \lambda^{\subs{V}}, \VM x^{\subs{V}} \right)$ meets a certain convergence criterium
\EndFor
\State $\hat \lambda = \lambda^{\subs{V}}$ and Ritz vector $\VM y = \VM V \VM x^{\subs{V}}$
\end{algorithmic}
\caption{ Arnoldi algorithm}  \label{algo:Arn}
\end{algorithm}

It is not difficult to prove that, in exact arithmetic, the eigenvalue approximations do not depend on the used continuation vector as long as the last element in $\VM t$ is non-zero. This changes when an error (an intended error or a numerical error) is introduced in one of the iterations. In these cases the choice of the continuation vector may be decisive for convergence. 

In the original version of the Arnoldi algorithm \cite{Arnoldi1951}, the $k$th standard basis vector is used as continuation vector in the $k$th iteration, but this choice performs badly if there is a significant error in one of the iterations, see \cite{Lee2007}. In the next section, another continuation vector is introduced which appears to be more robust against introduced errors. We motivate our choice for the continuation vector and give some convergence results. The convergence behaviour under errors is closely related to inexact Krylov methods, see \cite{Simoncini2003, Simoncini2005}. 
 
% \subsection{Continuation vector and Residual Arnoldi} \label{sect:resid_Arn}

In this section we want to answer the question which vector, after multiplying by $\VM A$, is the best to add to the subspace when this is a Krylov space contaminated with errors. This was the research topic in the papers \cite{Ye2008, Wu2014}. In both papers, the authors try to search for the continuation vector $\VM{t}$ that maximizes
\begin{equation} \label{eqn:VW_tk}
 \cos \left( \theta( \VM x, \text{Range} \{ \VM V, \VM A \VM V \VM t \} ) \right) := \max_{ \VM z \in \text{Range} \{ \VM V, \VM A \VM V \VM t \} } \cos( \theta( \VM x, \VM z))
\end{equation}
The following result is proved in \cite{Ye2008}. 
\begin{theorem} \label{the:optim_tk}
 Let $\VM R = \VM A \VM V - \VM V \left( \VM V^T \VM A \VM V \right)$, then the continuation vector $\VM t$ that maximizes \eqref{eqn:VW_tk} is 
\begin{equation} \label{eqn:optim_tk}
\VM t = \VM R^+ \VM x + \VM b
\end{equation}
with $\VM b$ in the nullspace of $\VM R$, $\VM R^+$ the pseudo-inverse of $\VM R$ and $\VM x$ the eigenvector associated with $\lambda$.
\end{theorem}

Although this result gives us insight in the choice of the optimal continuation vector, it is not computable as the formula contains the exact eigenvector. The author of \cite{Ye2008} argued that in practice, the eigenvector \begin{equation} \label{eqn:t_k}
\VM t = \VM x^{\subs{V}}.
\end{equation}
associated with the desired eigenvalue of the reduced problem $\VM V^T \VM A \VM V$ is a good practical choice for the continuation vector. In the remainder of this paper, this is used as continuation vector, despite the fact that in the case of large non-Hermitian matrices, it may be better to use the \emph{refined Ritz vector} \cite{Wu2014}. 

Before the suitability as continuation vector was discussed in \cite{Ye2008}, it was already observed that this choice leads to a more robust version of the Arnoldi algorithm and was called \emph{Residual Arnoldi}  \cite{Lee2007, LeeStewart2007}. Throughout this paper, we use the same name. The name comes from the property that the vector with which we expand the subspace is the residual vector as follows from Lemma \ref{lem:resid}. The following result is well-known in the literature \cite[Proposition 6.8]{Saad2011} and is frequently used in the analysis of Krylov methods although it holds for all subspace methods.

\begin{lemma} \label{lem:resid}
Let $\left( \lambda^\subs{V}, \VM y\right)$ be a Ritz vector of $\VM A$ with respect to the subspace $\subs{V}$, then the residual satisfies
\begin{equation} \label{eqn:resid}
\VM r = \left( \VM I - \VM V \VM V^T \right) \VM A \VM y
\end{equation}
so the residual is the result of applying matrix $\VM A$ on $\VM y$ orthogonalized with respect to $\VM V$.
\begin{eproof}
We can express
\begin{align}
\VM r & = \VM A \VM y - \lambda^{\subs{V}} \VM y \nonumber \\
& = \VM A \VM V \VM x^{\subs{V}} - \VM V \lambda^{\subs{V}} \VM x^{\subs{V}} \nonumber \\
& = \VM A \VM V \VM x^{\subs{V}} - \VM V \left( \VM V^T \VM A \VM V \right) \VM x^{\subs{V}} \nonumber \\
& = \left( \VM I - \VM V \VM V^T \right) \VM A \VM V \VM x^{\subs{V}}. \nonumber
\end{align}
\end{eproof}
\end{lemma}

Lemma \ref{lem:resid} implies that choosing \eqref{eqn:t_k} as continuation vector can simplify the Arnoldi algorithm (see Algorithm \ref{algo:Arn}) provided the norm of the residual is used as convergence criterium.  
% as in line \ref{line:res} in iteration $k$ we calculate the same vector as in line \ref{line:V_k} in iteration $k+1$. 
The Residual Arnoldi algorithm is shown in Algorithm \ref{algo:resid_Arn}.

\begin{algorithm}[h] 
\hspace*{\algorithmicindent} \textbf{Input:} Matrix $\VM{A}$, initial vector $\VM{u} \in \R^n$ \\
\hspace*{\algorithmicindent} \textbf{Output:} Approximate eigenpair $(\hat \lambda, \VM{y})$ of $\VM A$ 
\begin{algorithmic}[1]
\State $\VM r= \VM A \VM u - \lambda^{\VM u} \VM u$, with $\lambda^{\VM u}$ the Ritz value associated to the subspace .
\For{$k = 1$ to $n$}
\State add $\VM r$ to the subspace $\subs{V}$ \label{line:V_k_resarn}
\State Calculate eigenpair $\left( \lambda^{\subs{V}}, \VM x^{\subs{V}} \right)$ of $\VM V^T \VM A \VM V$
\State $\VM y = \VM V \VM x^\subs{V}$ \label{line:w_resarn}
\State Calculate the residual $\VM r$ of $\left( \lambda^\subs{V}, \VM y\right)$ via Lemma \ref{lem:resid} \label{line:resid}
\State Stop if $\norm{ \VM r}$ is small enough
\EndFor
\State $\hat \lambda = \lambda^{\subs{V}}$ and Ritz vector $\VM y = \VM V \VM x^{\subs{V}}$
\end{algorithmic}
\caption{ Residual Arnoldi (RA) Algorithm}  \label{algo:resid_Arn}
\end{algorithm}

\subsection{Convergence of the Residual Arnoldi algorithm} \label{sect:conv_res_arn}

In this section, we model the error in the RA algorithm, we discuss some convergence results and give tolerances for the different errors such that convergence is assured.
We model the possible error $\VM e$ as an error on $\VM r$, i.e on the vector that we want to add to the subspace. The error $\VM e$ consists of two sources: on the one hand, the error $\VM e_1 $ on $\VM x^\subs{V}$ (before the orthogonalization on line \ref{line:w_resarn} in the RA algorithm) and, on the other hand, the error $\VM e_2$ on $\VM r$ (after the orthogonalization on line \ref{line:resid}). Note that $\VM e_2 $ is already orthogonal to $ \subs{V}$, because otherwise the added vector is no longer orthogonal to $\subs{V}$. To resume, the vector that we add to the subspace can be written as
\begin{align} 
\hat{ \VM r} & = \left( \VM I - \VM V \VM V^T \right) \VM A \VM V \left( \VM x^\subs{V} + \VM e_1\right) + \VM e_2  \nonumber \\
& = \VM r + \underbrace{ \underbrace{\left( \VM I - \VM V \VM V^T \right) \VM A \VM V \VM e_1}_{\VM e_1^\bot} + \VM e_2}_{\VM e}. \label{eqn:e_w_bot}
\end{align}

The convergence of this method was analysed in \cite{Lee2007, LeeStewart2007}. The authors assume that in each iteration $k$ the error $\xkb e$  relative to the residual, is bounded by some $\eta$, in the sense
\begin{equation}
\dfrac{ \norm{ \xkb e}}{ \norm{ \xkb{r}}} \leq \eta. \label{eqn:rel_error}
\end{equation}
The main result on which the analysis relies, is that in iteration $k$ there always exists $\VM B_k$ such that $k$ iterations of the original Arnoldi method on the matrix 
$$ \hat{ \VM A}_k := \VM{A} + \VM B_k,$$
in the absence of errors results in the same subspace as performing $k$ iterations of RA on matrix $\VM A$ where in each iteration an error $\xkb{e}$ is made. In these papers, the authors make the assumption that the size of $\norm{ \VM B_k}$ is proportional to $\eta$, i.e. there exist $C_\eta$ such that 
$$\norm{ \VM B_k} \leq \eta C_\eta.
$$
The authors indicate that they cannot prove this assumption, although they have some empirical evidence that it is satisfied. Other assumptions were made on the separation of the spectrum of both $\VM A$ and $ \hat{ \VM A}$.

\begin{theorem}
Assume that the above conditions are satisfied, let $C_X$ be a constant depending on $C_\eta$ and the spectrum of $\VM A$.  If $2 \eta C_X < 1$, then
$$\norm{ \xkb y - \VM x} \leq  \dfrac{\kappa_k}{1-2 \eta C_X}$$
where $\kappa_k$ is a bound for the difference between the approximation and the searched eigenvector when no error is made. This value goes to zero as $k$ increases.
\end{theorem}

We refer to \cite{LeeStewart2007} for a proof and more details. In practice this bound is not (or hardly) computable, so a workable value of $\eta$ needs to be used in practice.

If the method converges, then the question arises how it depends on $\eta$. Empirically, the authors of \cite{Lee2007, LeeStewart2007} noticed that the speed of convergence for $\eta$ lower than or equal to $10^{-3}$ stays unchanged but that for $1 \gg \eta > 10^{-3}$ the convergence can slow down a bit. For $\eta$ around $1$, we do not have convergence, which can be explained by the fact that the error is as large as the residual.

% \todo[inline]{ in \cite{Lee2007} staat: Specifically, let $r_k$ be the exact residual, and $\tilde{r}_k$ be the one used in subspace expansion. The introduced error is defined as $f_k = \tilde{rk} rk$. Empirically, we have observed that $f_k$ can be as large as $10^{-3}$, and the candidate approximations can still converge to the target with rate similar to the rate without error.}

We now bound the norm of the errors $\VM e_1$ and $\VM e_2$ such that convergence is still preserved, i.e. such that condition \eqref{eqn:rel_error} is fulfilled with $\eta < 10^{-3}$. Let 
\begin{equation} \label{eqn:e_X-e_R}
\dfrac{ \norm{ \VM e_1}}{\norm{ \VM x^\subs{V}}} = \norm{ \VM e_1} \leq \eta_X \text{ and } \dfrac{ \norm{ \VM e_2}}{ \norm{ \VM r + \VM e_1^\bot}} \leq \eta_R
\end{equation}
be tolerances for the relative error, then we obtain starting from \eqref{eqn:e_w_bot}

\begin{align*}
    \dfrac{ \norm{\VM e} }{\norm{ \VM r}} & \leq \dfrac{\norm{ \VM e_1^\bot} }{\norm{ \VM r}} + \dfrac{\norm{ \VM e_2} }{\norm{ \VM r}} \\
    & = \dfrac{ \norm{ (\VM I - \VM V \VM V^T) \VM A \VM V \VM e_1}}{\norm{ \VM r}} + \dfrac{\norm{ \VM e_2} }{\norm{ \VM r}} \\
    & \leq \underbrace{ \dfrac{ \norm{\VM A} \norm{ \VM e_1}}{\norm{ \VM r}}}_{\leq \eta/2} + \underbrace{ \dfrac{ \norm{ \VM e_2}}{\norm{ \VM r}}}_{\leq \eta/2}. 
\end{align*}

The aim is now to bound both terms by $\eta/2$. Remark that the first term means that the error propagation due to the matrix-vector product needs to be limited. From \eqref{eqn:e_X-e_R} it follows that 
\begin{equation} \label{eqn:bound_eX}
\eta_X  \leq \dfrac{ \eta \norm{ \VM r}  }{2 \norm{\VM A}} 
\end{equation}
% As
% $$
%      \dfrac{ \norm{\VM A} \norm{ \VM e_1}}{ \norm{ \VM r}} \leq  \dfrac{ \norm{\VM A} \eta_X}{\norm{ \VM r}}
% $$
% is smaller than $\epsilon/2$,it follows that

which means that $\eta_X$ scales with the residual norm and the norm of the matrix $\VM A$. The lower the residual is, the lower $\eta_X$ needs to be. In the example in Section \ref{sect:Error_an}, we demonstrate that if we keep $\eta_X$ constant, then convergence stagnates if the residual becomes small.

Using $\norm{ \VM r + \VM e^\bot_1} \leq \norm{ \VM r} + \norm{ \VM e^\bot_1}$, we obtain
\begin{align*}
\dfrac{ \norm{ \VM e_2}}{\norm{ \VM r}} & = \dfrac{ \norm{ \VM e_2}}{ \norm{ \VM r + \VM e_1^\bot}}  \dfrac{ \norm{ \VM r + \VM e_1^\bot} }{\norm{ \VM r}}  \\
    & \leq \dfrac{ \norm{ \VM e_2}}{ \norm{ \VM r + \VM e_1^\bot}} \dfrac{ \norm{ \VM r} + \norm{e^\bot_1} }{\norm{ \VM r}} \\
    & \leq \dfrac{ \norm{ \VM e_2}}{ \norm{ \VM r + \VM e_1^\bot}} (1 + \eta/2).
\end{align*}
As this needs to be smaller than $\eta/2$, we must impose that
\begin{equation} \label{eqn:bound_eR}
    \eta_R \leq  \dfrac{\eta }{2 + \eta}.
\end{equation}

%Using $\norm{ \VM r + \VM e^\bot_1} \leq \norm{ \VM r} + \epsilon/2$ we obtain 
%\begin{align*}
%    \dfrac{ \norm{ \VM e_2}}{\norm{ \VM r}} & = \dfrac{ \norm{ \VM e_2}}{ \norm{ \VM r + \VM e_1^\bot}}  \dfrac{ \norm{ \VM r + \VM e_1^\bot} }{\norm{ \VM r}}  \\
%    & \leq \dfrac{ \norm{ \VM e_2}}{ \norm{ \VM r + \VM e_1^\bot}} \dfrac{ \norm{ \VM r} + \epsilon/2}{\norm{ \VM r}} \\
%    & \leq \eta_R \dfrac{ \norm{ \VM r} + \epsilon/2}{\norm{ \VM r}}.
%\end{align*}
%As this term needs to be smaller than $\epsilon/2$, we must impose that
%\begin{equation} \label{eqn:bound_eR}
%    \eta_R \leq  \dfrac{\epsilon \norm{ \VM r}}{2 \left( \norm{ \VM r} + \epsilon/2\right)} \leq \epsilon/2.
%\end{equation} 
%\textcolor{red}{
%---------------------------------------------------------------------------------------------------------------------------------------------------------------------
%Hetgeen hierboven staat klopt niet. Using $\norm{ \VM r + \VM e^\bot_1} \leq \norm{ \VM r} + \norm{ \VM e^\bot_1}$, we obtain
%\begin{align*}
%\dfrac{ \norm{ \VM e_2}}{\norm{ \VM r}} & = \dfrac{ \norm{ \VM e_2}}{ \norm{ \VM r + \VM e_1^\bot}}  \dfrac{ \norm{ \VM r + \VM e_1^\bot} }{\norm{ \VM r}}  \\
%    & \leq \dfrac{ \norm{ \VM e_2}}{ \norm{ \VM r + \VM e_1^\bot}} \dfrac{ \norm{ \VM r} + \norm{e^\bot_1} }{\norm{ \VM r}} \\
%    & \leq \dfrac{ \norm{ \VM e_2}}{ \norm{ \VM r + \VM e_1^\bot}} (1 + \eta/2).
%\end{align*}}
%\textcolor{red}{
%
%---------------------------------------------------------------------------------------------------------------------------------------------------------------------} 

\section{Tensor-Arnoldi for one parameter} \label{sec:Tensor-Arn_1}
We are ready to propose our algorithm for calculating the wanted eigenvalue over the parameter space. First, a general description of the algorithm is given, then we describe in more detail some specific steps of the algorithm. We explain the algorithm for the case where we have one parameter. For the extension to multiple parameters we refer to Section \ref{sect:Tensors}. First, we discretize $\Omega$ into $n_1$ points $\omega_1, \omega_2, \hdots, \omega_{n_1}$. We want to compute the desired eigenvalue $\lambda$ at these points. For this, we consider the linear operator $\Ten A$ such that

% \todo[inline]{zeggen dat we de basis willen reeel houden, hoewel $\VM X$ mogelijks complex is}
\begin{equation} \label{eqn:diagAX}
\Ten A \VM x = \begin{bmatrix}
\VM A(\omega_1) & & \\
& \ddots & \\
& & \VM A(\omega_{n_1})
\end{bmatrix} \begin{bmatrix}
\VM x(\omega_1) \\
\VM x(\omega_2) \\
\vdots \\
\VM x(\omega_{n_1})
\end{bmatrix} = \begin{bmatrix}
\VM A(\omega_1) \VM x(\omega_1) \\
\vdots \\
\VM A(\omega_{n_1}) \VM x(\omega_{n_1})
\end{bmatrix} \in \R^{n n_1}.
\end{equation}
If $\VM x$ is reformulated as matrix $\VM X = [\VM x(\omega_1), \VM x(\omega_2), \hdots, \VM x(\omega_{n_1})] \in \C^{n \times n_1}$, then \eqref{eqn:diagAX} is equivalent to an operator $\Ten{A}( \VM X)$ such that $\text{vec}( \Ten{A}( \VM X)) = \Ten A \VM x$. Conceptually, our method boils down to applying RA to \eqref{eqn:diagAX} which is nothing more than applying RA for each $\omega_i$ individually.

In order to efficiently compute $\Ten A(\VM X)$ we need the expansion in \eqref{eqn:VW_A}. If we have a decomposition $\VM X = \VM U \VM Z$ with $\VM U \in \R^{n \times n_U}$ and $\VM Z \in \C^{n_U \times n_1}$, then in the same way as in \cite{Kressner2011} we obtain

\begin{equation} \label{eqn:A(X)}
\Ten{A}(\VM X) = \underbrace{[\VM A_1 \VM U, \VM A_2 \VM U, \hdots, \VM A_m \VM U]}_{\VM W_U \in \R^{n \times n_U m}} \underbrace{\begin{bmatrix}
\VM Z \cdot \text{diag}( [f_1(\omega_1), f_1(\omega_2), \hdots, f_1(\omega_{n_1})]) \\
\VM Z \cdot \text{diag}( [f_2(\omega_1), f_2(\omega_2), \hdots, f_2(\omega_{n_1})])\\
\vdots \\
\VM Z \cdot  \text{diag}( [f_m(\omega_1), f_m(\omega_2), \hdots, f_m(\omega_{n_1})])\\
\end{bmatrix}}_{\VM W_Z \in \C^{n_U m \times n_1}}.
\end{equation}

Remark that the term $$\VM Z \cdot \text{diag}( [f_i(\omega_1), f_i(\omega_2), \hdots, f_i(\omega_{n_1})]) = [ f_i(\omega_1) \VM Z(:,1), f_i(\omega_2) \VM Z(:,2), \hdots, f_i(\omega_{n_1}) \VM Z(:,n_1)],$$ so column $j$ of $\VM Z$ is multiplied with the function $f_i$ evaluated in $\omega_j$.
Let now $\VM X^\subs{V}$ be the matrix storing in each column $i$ the eigenvector $\VM x^\subs{V}(\omega_i)$ of the reduced eigenvalue problem $\VM V^T \VM A(\omega_i) \VM V$ with respect to the subspace $\subs{V}$. We first make a low-rank decomposition $\VM U \VM Z$ of $\VM X^\subs{V}$ and subsequently we compute $\VM W_U \VM W_Z= \Ten{A}( (\VM V \VM U) \VM Z)$ via \eqref{eqn:A(X)} and we orthogonalise the result with respect to $\subs{V}$ to obtain an approximate residual for each parameter value, see Lemma \ref{lem:resid}. This is analogeous to lines \ref{line:w_resarn} and \ref{line:resid} in the RA algorithm (see Algorithm \ref{algo:resid_Arn}). Instead of adding all column vectors of the residual, we make a low-rank approximation and add in the next iteration the resulting basis $\VM U$ to $\subs{V}$. The algorithm for the low-rank decomposition is the same as in \cite{Kressner2011}. 

% \todo[inline]{De keuze voor $\epsilon$ als tolerantie is ongelukkig, de fout die we maken hebben we immers ook $\epsilon$ genoemd. In thesis heb ik de tolerantie voor de fout hernoemd naar $\eta$}
After calculating new eigenvalues and eigenvectors we can repeat the same procedure untill all residuals are below a certain tolerance $\epsilon$. As we are adding multiple vectors to the subspace in each iteration and the maximal dimension of the subspace is assumed to be prescribed, the question arises how we restart the subspace when it attains its maximal dimension. In the literature, several techniques are developed for the Arnoldi algorithm like polynomial restart methods \cite{Saad1980, Saad1984} or the implicitly restarted method \cite{Sorensen1992}. All of them try to delete the part of the subspace that corresponds to an unwanted part of the spectrum, but still keeps a significant part of the subspace. The reason is that not only the Ritz vector associated with the wanted eigenvalue, but also some other Ritz vectors whose eigenvalues are close to $\lambda^{\subs{V}}$ help improving the convergence to the wanted eigenvalue. This is known as \emph{thick restarting}. In our case of parametric eigenvalue problems, it is slightly more complicated as eigenvectors vary over the parameter. This means that for just representing the wanted eigenvalue we already need quite a large subspace. If also eigenvectors associated with other eigenvalues need to be present in the restarted subspace, there is less room for iterations within the same subspace. Therefore, we choose to just use a basis for the first eigenvectors as restarted subspace. This is a good strategy since it is known from \cite{ruymbeek2019calculating, seyranian2003multiparameter} that the derivative of an eigenvector, which is a measure for the change of the eigenvector, depends the most on the eigenvectors whose eigenvalues are the closest to $\lambda(\VM \omega)$. As we estimate the eigenvector over the parameter space, we therefore also estimate the eigenvector whose eigenvalues are the closest to $\lambda(\VM \omega)$. The algorithm can be found in Algorithm \ref{algo:resid_Arn_param}.
% Moreover we are in fact adding more information than just one Ritz vector as it is known that the change in the eigenvector over the parameter is affected mostly by the eigenvectors whose eigenvalues are close to the wanted eigenvalue, see also \cite{ruymbeek2019calculating, seyranian2003multiparameter}. 

\begin{algorithm}[h] 
\hspace*{\algorithmicindent} \textbf{Input:} $\VM{A}(\omega)$, orthonormal matrix $\VM U \in \R^{n \times p}$, tolerance $\epsilon$, tolerances $\eta_X$ and $\eta_R$ for the two decompositions, maximal dimension of subspace $n^\subs{V}$ \\
\hspace*{\algorithmicindent} \textbf{Output:} Approximate eigenpair $(\hat \lambda(\omega_i), \VM{y}(\omega_i)), i = 1,2, \hdots, n_1$
\begin{algorithmic}[1]
\State $\hat{ \VM R} = [1, 1, \hdots, 1]$, \quad $k = 1$
\While{$\norm{ \hat{ \VM R}}_F/\sqrt{n_1} > \epsilon$} 
\State Add $\VM U$ to the subspace $\subs{V}$
\State Calculate eigenpair $\left( \lambda^{\subs{V}}(\omega_i), \VM x^{\subs{V}}(\omega_i) \right)$ of $\VM V^T \VM A (\omega_i) \VM V, i = 1,2, \hdots, n_1$ \label{line:reduc_eigv_probl}
\State $\VM X^{\subs{V}} = [\VM x^{\subs{V}}(\omega_1), \VM x^{\subs{V}}(\omega_2), \hdots, \VM x^{\subs{V}}(\omega_{n_1})]$  \label{line:XV}
\If{ dimension $\subs{V} > n^\subs{V}$}
\State $[\VM U, \VM Z] = \text{lowrank}( \VM X^{\subs{V}}, 10^{-12})$ \label{line:lowrank1a} 
\State Make $\subs{V}$ empty
\Else
\State $[\VM U, \VM Z] = \text{lowrank}( \VM{ X}^{\subs{V}}, \eta_X)$ \label{line:lowrank1b} 
% \State $ \hat{ \VM{X}}^{\subs{V}_k} = \xkb U \xkb Z$
\State $[\VM W_U, \VM W_Z] = \Ten{A}( \left( \VM V \VM U \right) \VM Z)$ using \eqref{eqn:A(X)} \label{line:A(X)}
\State $\VM W^\bot_U = (\VM I - \VM V \VM V^T) \VM W_U$ ( Orthogonalize $\VM W_U$ w.r.t $\VM V$)
\State Residual $ \hat{ \VM R} := [\hat{\VM{r}}(\omega_1), \hat{\VM{r}}(\omega_2), \hdots, \hat{\VM{r}}(\omega_{n_1})] = \VM W^\bot_U \VM W_Z$ \label{line:residual}
\State $[\VM U, \sim] = \text{lowrank}( \VM W^\bot_U, \VM W_Z, \eta_R)$; \label{line:lowrank2}
\EndIf
\State $k = k+1$ 
\EndWhile
\State $\hat \lambda(\omega_i) = \lambda^{\subs{V}}(\omega_i)$ and Ritz vector $\VM y(\omega_i) = \VM V \VM X^{\subs{V}}(\omega_i), i= 1,2,\hdots,n_1$
\end{algorithmic}
\caption{ Parametric Residual Arnoldi (PRA) Algorithm}  \label{algo:resid_Arn_param}
\end{algorithm}

For the calculation of $\VM V^T \VM A(\omega_i) \VM V, i = 1, \hdots, n_1$, we first calculate $\VM V^T \VM A_i \VM V, i = 1, \hdots, m$ and sum these smaller matrices appropriately. The operation $\VM A_i \VM V, i = 1,\hdots m$ is not expensive as these $\VM A_i$ is usually sparse, but the inner product $\VM V^T \left( \VM A_i \VM V\right)$ need $2 n^2 n_v^2$ operations with $n_v$ the dimension of $\VM V$. As $\VM V$ is built iteratively, results from the previous iteration can be used if the subspace was not restarted. The lower $\eta_X$ and $\eta_R$, the higher the number of these inner products. In the example in the Section \ref{sect:Error_an}, we count the number of these inner products and compare them for different values of $\eta_X$ and $\eta_R$.

\section{Error analysis} \label{sect:Error_an}

We illustrate the performance of the proposed algorithm through an example. We show here what the impact is of the choice of the tolerance levels $\eta_R$ resp. $\eta_X$ in the PRA algorithm and of the maximal size of the subspace $\subs{V}$. In particular, we explain why choosing $\eta_R$ relatively high gives the best results, i.e we need less iterations and less expensive matrix-vector products without affecting the convergence.

\begin{figure}[h]
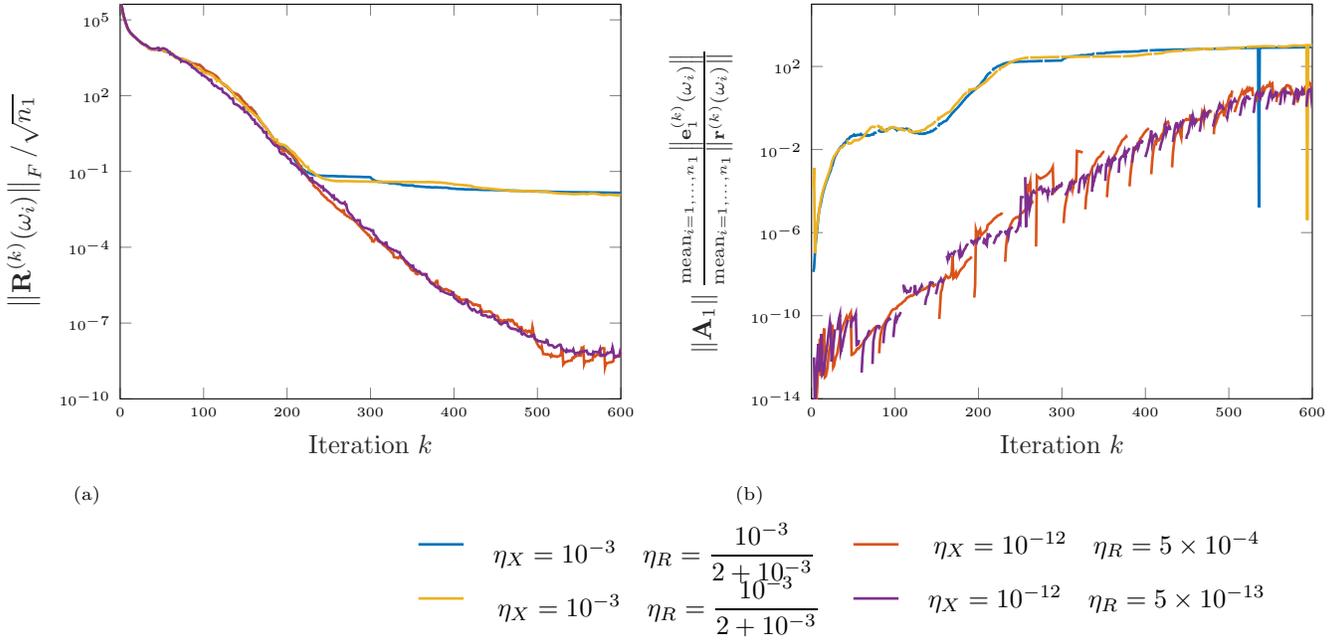

\begin{center}
\setlength{\figW}{7cm} % determines the width of the images
  \hspace{-0.5\linewidth}
 \begin{subfigure}[h]{0.3\textwidth}
%  \centering
  \hspace{1.6cm}
  \input{Afbeeldingen/rel_residu_conv_diff.tex}
  \caption{}
\end{subfigure} 
\hspace{0.15\linewidth}
 \begin{subfigure}[h]{0.3\textwidth}
%  \centering
  \hspace{1.6cm}
  \input{Afbeeldingen/rel_E1_conv_diff.tex}
  \caption{}
\end{subfigure} \\
\begin{tikzpicture}
\hspace{0.05\linewidth}
\begin{customlegend}[legend columns=2,legend style={align=left,draw=none,column sep=2ex},legend entries={$\eta_X= 10^{-3} \quad \eta_R = \dfrac{ 10^{-3} }{2 + 10^{-3} }$ , $\eta_X= 10^{-12} \quad \eta_R = 5 \times 10^{-4}$, $\eta_X= 10^{-3} \quad \eta_R = \dfrac{ 10^{-3} }{2 + 10^{-3} }$, $\eta_X= 10^{-12} \quad \eta_R = 5 \times 10^{-13}$}]
\addlegendimage{color=mycolor1, line width=1.0pt,line legend}
\addlegendimage{color=mycolor2,solid,line width=1.0pt}  
\addlegendimage{color=mycolor3,solid,line width=1.0pt}
\addlegendimage{color=mycolor4,solid,line width=1.0pt}
\end{customlegend}
\end{tikzpicture}
\end{center}
\caption{Results of the PRA algorithm for different tolerances $\eta_X$ and $\eta_R$ with maximal subspace dimension $150$ for Example \ref{ex:conv_diff}.  In (a) we plot the Frobenius norm of the residual, scaled with the number of sample points. In (b) the average first error $\VM e_1$ relative to the average residual, scaled with the norm of $\VM A_0$ is shown. As long as the value is below or around $10^{-3}$, the algorithm converges as if there were no error, see the criterium in \eqref{eqn:bound_eX}. If the value is large $(> 1)$ convergence is lost. } \label{fig:conv_diff_results}
\end{figure}

\begin{eexample} \label{ex:conv_diff}
Consider the convection-diffusion differential equation
\begin{equation} \label{eqn:conv_diff}
f = c_1 u_x + c_2 u_y + d_{11} u_{xx} + 2d_{12} u_{xy} + d_{22} u_{yy} 
\end{equation}
for $(x, y) \in \Xi= (0, 1) \times (0, 1)$. The real function $u$ has the argument $(x, y)$. 
The constants $c_1$ and $c_2$ are associated with the convection part and the constants $d_{11}$, $d_{12}$ and $d_{22}$ are associated with the diffusion part of the equation. It is always assumed that matrix
$$\VM D = \begin{bmatrix}
d_{11} & d_{12} \\
d_{12} & d_{22} \end{bmatrix}
$$
is positive semi-definite. For the discretisation of $\Xi$ we use second-order central finite difference formulas and we assume Dirichlet boundary conditions. In this way, we obtain a matrix $\VM A$ depending on the values $c_1$,$c_2$,$d_{11}$,$d_{12}$, $d_{22}$ such that 
\begin{equation} \label{eqn:conv_diff_A}
\VM A(c_1, c_2, d_{11}, d_{22}, d_{12}) = c_1 \VM{ \hat{ A}}_1 + c_2 \VM{ \hat{ A}}_2 + d_{11} \VM{ \hat{ A}}_3 + d_{22} \VM{ \hat{ A}}_4 + d_{12} \VM{ \hat{ A}}_5.
\end{equation} 
For this first example we take $c_2 = 1, d_{11} = 1.1, d_{22} = 1, d_{12} = 1$ and we let $c_1$ vary in the interval $\Omega = [-2.5, 2.5]$. We discretize this interval in $100$ points. This means we can rewrite \eqref{eqn:conv_diff} as 
$$\VM A(c_1) = \VM A_1 + c_1 \VM A_2.$$
 The desired eigenvalue is the eigenvalue with largest real part, in this case the same as the eigenvalue closest to zero. The region $\Xi$ is discretised into $100$ points in each dimension, so, the resulting matrices $\VM A_i$ are of size $10^4 \times 10^4$. The results can be found in Figure \ref{fig:conv_diff_results} and Table \ref{tbl:conv_diff}. Remark that we need many iterations but tests with individual parameters values reveal that also the Arnoldi algorithm needs many iterations.
\end{eexample}
In Section \ref{sect:conv_res_arn} we observed that, in the RA algorithm, we add the residual in each iteration to the subspace and conclude that this is more robust than classical Arnoldi. In the PRA algorithm, however, the subspace is enriched with an approximate residual for all sample points by making low-rank approximations. For this, we rely on the error analysis made in Section \ref{sect:conv_res_arn}. At line~\ref{line:lowrank1b}, on the one hand, and at line \ref{line:lowrank2} of PRA algorithm on the other, we make low-rank approximations of the matrix consisting of the eigenvector resp. the approximate residual and so introducing errors $\VM E_1 := [\VM e_1(\omega_1), \VM e_1(\omega_2), \hdots, \VM e_1(\omega_{n_1})]$ and $\VM E_2 := [\VM e_2(\omega_1), \VM e_2(\omega_2), \hdots, \VM e_2(\omega_{n_1})]$.
We first discuss the error $\VM E_1$. The low-rank approximation is such that the relative error is bounded by the tolerance $\eta_X$, i.e
\begin{equation} \label{eqn:eps_X}
\dfrac{ \norm{ \VM E_1}_F}{\norm{ \VM X^\subs{V}}_F} \leq \eta_X.
\end{equation}
Since $\norm{ \VM X^\subs{V}}_F = \sqrt{n_1}$, it follows that $\norm{ \VM E_1}_F \leq \eta_X \sqrt{n_1}$. 
Technically, we can only prove from this that 
$$\norm{ \VM e_1(\omega_i)}^2 \leq \eta_X^2 n_1, i = 1, 2, \hdots, n_1,$$
but assuming the error is uniformly distributed over the sample points, we can assume that
\begin{equation} \label{eqn:e_X}
\norm{ \VM e_1(\omega_i)}  \lessapprox \eta_X.
\end{equation}

Figure \ref{fig:conv_diff_results} (a) shows the average residual over the parameter and we observe that for $\eta_X = 10^{-3}$ the convergence slows down after $400$ iterations. This is related to the average relative error $\VM e_1$ plotted in Figure \ref{fig:conv_diff_results} (b), as this becomes too high such that condition \eqref{eqn:bound_eX} is not fulfilled.  
Analoguously to \eqref{eqn:e_X} we can prove that $$\dfrac{ \norm{ \VM e_2(\omega_i)}}{\norm{ \VM r(\omega_i) + \VM e_1^\bot }} \lessapprox \eta_R,  i = 1, 2, \hdots, n_1.$$
As seen in \eqref{eqn:bound_eR}, we only need to make sure that $\eta_R \leq \dfrac{\eta}{2 + \eta}$, where $\eta$ can be chosen as high as $10^{-3}$. This explains why in Figure \ref{fig:conv_diff_results} (a) choosing $\eta_R = \dfrac{10^{-3}}{2+10^{-3}}$ does not harm convergence. Furthermore, it seems even that convergence is faster for higher $\eta_R$.
This phenomenon is intriguing and cannot be explained only by a relative error plot. It is generally known in Arnoldi-type algorithms that, with each restart, convergence can slow down. The higher the tolerance $\eta_R$ the less vectors are added to the subspace, so we do not need to restart so often, hence we can do more iterations within the same subspace and convergence is faster. 

In Table \ref{tbl:conv_diff}, we see statistics about the computational cost. As expected, we observe that the number of restarts increases as the tolerances decrease and that the average size of the subspace is independent of the chosen tolerances. As indication of the increased computational cost when taking the tolerances smaller, we count the number of inner products when creating the reduced eigenvalue problem (line \ref{line:reduc_eigv_probl} in the PRA algorithm). A more detailed explanation on this was given at the end of Section \ref{sec:Tensor-Arn_1}. 
\begin{table}[h]
\begin{center}
\begin{tabular}{|l|l|l|l|l|}
\hline
&$\epsilon_X=10^{-3}$ &$\epsilon_X=10^{-12}$ &$\epsilon_X=10^{-3}$ &$\epsilon_X=10^{-12}$ \\
& $\epsilon_R =5 \times 10^{-4}$ & $\epsilon_R =5 \times 10^{-4}$ & $\epsilon_R =5 \times 10^{-13}$ & $\epsilon_R =5 \times 10^{-13}$ \\\hline
Number of restarts&$7$&$16$&$15$&$38$\\\hline
Avg. size of subspace&$79.5233$&$78.7383$&$79.545$&$81.0283$\\\hline
Avg. number of inner products &$691.74$&$1429.9167$&$1272.9133$&$3434.36$\\
for reduced eigenvalue problem & & & & \\\hline
\end{tabular}
\end{center}
\caption{Statistics associated to the results shown in Figure \ref{fig:conv_diff_results}.}
% This file was created in matlab file table_conv_diff on 18-Feb-2020
 \label{tbl:conv_diff}
\end{table} 

\begin{figure}[h]
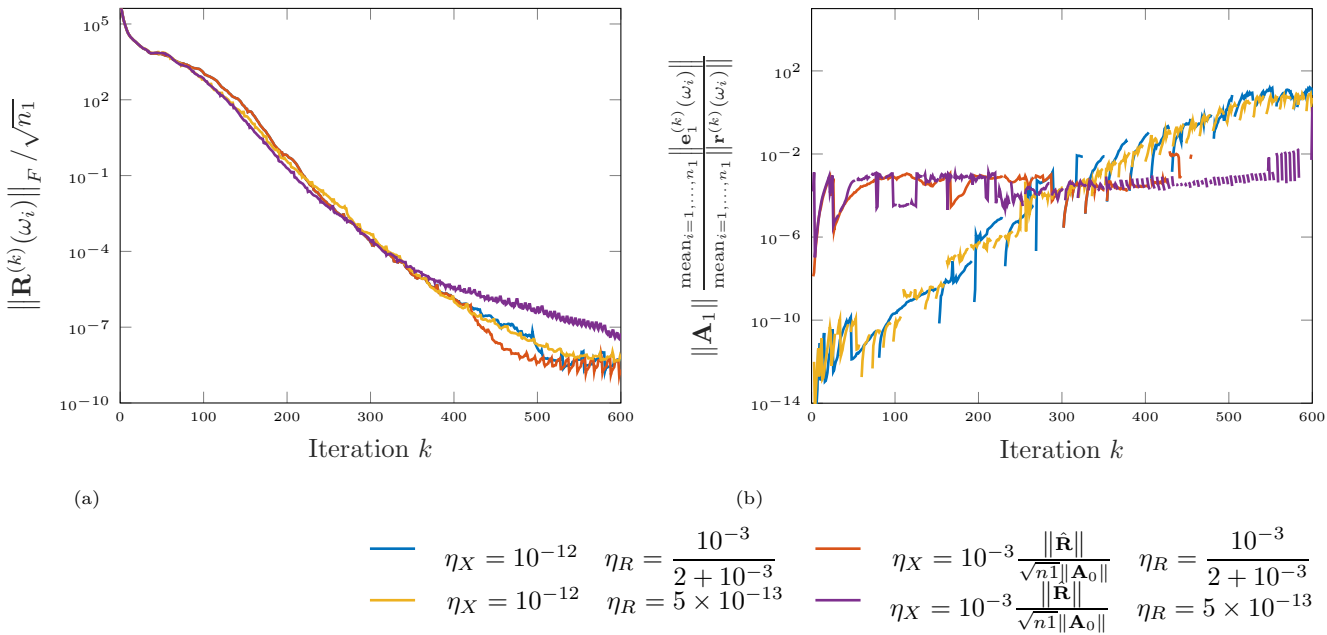

%\end{subfigure} 
\begin{center}
\setlength{\figW}{7cm} % determines the width of the images
  \hspace{-0.5\linewidth}
 \begin{subfigure}[h]{0.3\textwidth}
%  \centering
  \hspace{1.6cm}
  \input{Afbeeldingen/rel_residu_conv_diff_vartol_X.tex}
  \caption{}
\end{subfigure} 
\hspace{0.15\linewidth}
 \begin{subfigure}[h]{0.3\textwidth}
%  \centering
  \hspace{1.6cm}
  \input{Afbeeldingen/rel_E1_conv_diff_vartol_X.tex}
  \caption{}
\end{subfigure} \\
\begin{tikzpicture}
\hspace{0.05\linewidth}
\begin{customlegend}[legend columns=2,legend style={align=left,draw=none,column sep=2ex},legend entries={$\eta_X= 10^{-12} \quad \eta_R = \dfrac{10^{-3}}{2+10^{-3}}$ , $\eta_X= 10^{-3} \frac{ \norm{\hat{ \VM R}} }{\sqrt{n1} \norm{ \VM A_0}}
 \quad \eta_R = \dfrac{10^{-3}}{2+10^{-3}}$, $\eta_X= 10^{-12} \quad \eta_R = 5 \times 10^{-13}$, $\eta_X= 10^{-3} \frac{ \norm{\hat{ \VM R}} }{\sqrt{n1} \norm{ \VM A_0}}
 \quad \eta_R = 5 \times 10^{-13}$}]
\addlegendimage{color=mycolor1, line width=1.0pt,line legend}
\addlegendimage{color=mycolor2,solid,line width=1.0pt}  
\addlegendimage{color=mycolor3,solid,line width=1.0pt}
\addlegendimage{color=mycolor4,solid,line width=1.0pt}
\end{customlegend}
\end{tikzpicture}
\end{center}
\caption{Results of the PRA algorithm for different tolerances $\eta_X$ and $\eta_R$ with maximal subspace dimension $150$ for Example \ref{ex:conv_diff}. We compare the results when taking a constant low tolerance $\eta_X$ with the case that $\eta_X$ is relative to the residual. In (a) we plot the Frobenius norm of the residual, scaled with the number of sample. We see that in all four cases the algorithm converges. In (b), the average error $\VM e_1$ relative to the average residual, scaled with the norm of $\VM A_1$ is shown. As expected, we see that this error relative to the residual stays constant. } \label{fig:conv_diff_results_vartol}
\end{figure}

In the experiments in Figure \ref{fig:conv_diff_results}, we observed that if $\eta_X$ is chosen high, this harms the convergence when the norm of the error $\VM e_1$ starts to be high relative to the residual. The solution for this is to choose the tolerance $\eta_X$ in iteration $k$ relative to the residual in the previous iteration. In Figure \ref{fig:conv_diff_results_vartol}, we compare results when holding $\eta_X$ constant and when we adapt its value relatively with respect to the residual. We see that we obtain similar results. This is explained by Figure \ref{fig:conv_diff_results_vartol}(b) as we see that the error $\VM e_1$ relative to the residual remains constant.

The conclusion is that we obtain the best results when taking $\eta_R$ high and taking $\eta_X$ relative to the residual. In this way the norm of both errors $\VM e_1$ and $\VM e_2$ stay constant with respect to the residual, which preserves convergence if this constant is low enough ($< 10^{-3}$), see Section \ref{sect:Error_an}. 

\section{More than one parameter} \label{sect:Tensors}

In this section, we extend the PRA algorithm to the case of more than one parameter. Not only all computations need to have an equivalent tensor counterpart, we also need a suitable representation of the tensors as the storage requirement grows exponentially with the dimension. In this section, we first introduce some notation, then we discuss the chosen tensor-format and explain how all operations are adapted. 

We assume that $\Omega$ is a Cartesian product of intervals in $\R$, so $$\Omega = \Omega_1 \times \Omega_2 \times \hdots \times \Omega_d = [a^1, b^1] \times [a^2, b^2] \times \hdots \times [a^d, b^d].$$ Each dimension $i, i=1, 2, \hdots, d$ is discretised into $n_i$ points, making a total of $n_0 = n_1 n_2 \hdots n_d$ sample points where we want to approximate the desired eigenvalue. A point $\VM \omega$ on this grid is denoted as $\left( \omega_{i_1}, \omega_{i_2}, \hdots, \omega_{i_d} \right)$ for $i_j = 1,2, \hdots, n_j, j = 1,2, \hdots, d$. 
All tensor equivalents are the calligraphic version of the letters used in the one parameter case. For example, the eigenvectors are stored in tensor $\Ten X \in \C^{n \times n_1 \times \hdots \times n_d}$ where 
$$\Ten X(:,i_1, i_2, \hdots, i_d) = \VM x(\omega_{i_1}, \omega_{i_2}, \hdots, \omega_{i_d}).$$
We define tensor $\Ten{F} \in \R^{m \times n_1 \times n_2 \times \hdots \times n_d}$ as the tensor with all function values, so  
$$\Ten{F}(i_0,i_1, \hdots, i_d) = f_{i_0} \left( \omega_{i_1}, \omega_{i_2}, \hdots, \omega_{i_d} \right)$$
for a point $\VM \omega =  \left( \omega_{i_1}, \omega_{i_2}, \hdots, \omega_{i_d} \right).$ If $d = 1$, this is the matrix 
$$\VM F = \begin{bmatrix}
f_1(\omega_1) & f_1(\omega_2) & \hdots & f_1(\omega_{n_1}) \\ 
f_2(\omega_1) & f_2(\omega_2) & \hdots & f_2(\omega_{n_1}) \\ 
\vdots \\
f_m(\omega_1) & f_m(\omega_2) & \hdots & f_m(\omega_{n_1}) \\ 
\end{bmatrix}.$$

\subsection{Tensors}
We first state some basic operations that can be performed on tensors, then we describe how to adapt the PRA Algorithm. A detailed overview on tensors can be found in \cite{Kolda2009}.

\begin{definition}
\begin{enumerate}
\item The \textit{mode-$i$ matriciziation} or \textit{unfolding} of a tensor $\Ten{Z} \in \C^{n \times n_1 \times \hdots \times n_d}$ is a matrix $\VM{Z}_{(i)} \in \C^{n_i \times (n n_1 \hdots n_{i-1} n_{i+1}\hdots n_d)}$ where the columns are the mode-$i$ vectors. The order in which they are placed in  $\VM{Z}_{(i)}$ is chosen to be canonical.
\item The \textit{mode-$i$ multiplication} of a matrix $\VM{M} \in \R^{k \times n_i}$ with tensor $\Ten{Z}$ is defined as 
$$\Ten{Y} = \left( \VM{I},  \hdots, \VM{I},\VM{M},\VM{I}, \hdots, \VM{I} \right) . \Ten{Z} \Leftrightarrow \VM{Y}_{(i)} = \VM{M} \VM{Z}_{(i)}$$  
where $\VM{M}$ is at the $i+1$th position in the tuple.
\item The \emph{Hadamard product} between two tensors $\Ten U$ and $\Ten W$ of the same dimension, is obtained by the elementwise multiplication and is denoted by $\Ten U ._H \Ten W$.
\end{enumerate}
\end{definition}

If $\VM U$ is the column basis for the tensor $\Ten X$, then we can write
\begin{equation} \label{eqn:decomp}
\Ten{X} := \left( \VM{U}, \VM{I}, \VM{I}, \hdots, \VM{I} \right) . \Ten{Z}
\end{equation} 
with $\VM{U} \in \R^{n \times k}$ and $\Ten{Z} = \R^{k \times n_1 \times n_2 \times \hdots \times n_d}$, where $\VM{U}$ is a basis for $\VM{X}_{(0)}$.

We look now for a suitable tensor representation such that all the computations in the iteration can be done in an efficient way. A strong requirement for the decomposition of $\Ten{X}^\subs{V}$ besides the storage requirement is that the following operations in PRA Algorithm can be computed efficiently within this decomposition:
\begin{enumerate}
\item Calculating a real basis for the column space, even if the tensor contains complex numbers (line \ref{line:lowrank1a}, \ref{line:lowrank1b} and \ref{line:lowrank2} in the PRA Algorithm). \label{enum:2}
\item Approximating a tensor up to a given tolerance (line \ref{line:lowrank1a}, \ref{line:lowrank1b} and \ref{line:lowrank2}), in the same sense as in \eqref{eqn:eps_X}.
\item Computing a decomposition of $\Ten{W}_Z$ given the decomposition of $\Ten{Z}$ (line \ref{line:A(X)}). \label{enum:1}
\end{enumerate} 

% \todo[inline]{ $\R$ moet hier $\C$}
The Tensor-Train format (TT-format) satisfies all the above conditions. Before we illustrate this we briefly recall the TT-format. A more detailed overview can be found in \cite{Oseledets2011}.
A tensor $\Ten{Y} \in \C^{n_0 \times n_1 \times \hdots \times n_d}$ is in TT-format if  
\begin{equation}
\Ten{Y}(i_0,i_1,\hdots,i_d) = \VM G^Y_0(i_0) \VM G^Y_1(i_1) \hdots \VM G^Y_d(i_d)    
\end{equation}
with $\VM G^Y_j(i_j)$ a $r_{j-1}^Y \times r_j^Y$ matrix. 
The length $r^Y_{j-1}$ and width $r^Y_{j}$ of these matrices $r^Y_i, i= -1,0,\hdots, d$ are called the \emph{TT-ranks}. It can immediately be observed that $r^Y_{-1} = r^Y_{d} = 1$. We call $r^Y_i, i = 0,1, \hdots, d-1$ the \emph{interior} ranks.
The storage of the TT-representation is lower than $(d+1)n r^2$ with $r^Y_i \leq r$ and $n_i \leq n$, which is for $r$ small significantly smaller than the storage of the full tensor. The matrices $\VM G^Y_i$ can be stored in a tensor $\Ten{G}^Y_j \in \R^{r^Y_{j-1} \times n_j \times r^Y_j}, j=0,1,\hdots,d$ as $\Ten{G}^Y_j(:,i_j,:) = \VM G^Y(i_j), j = 1,2, \hdots, n_j$.

For computing a TT-decomposition of a tensor, we perform for each mode consecutively a matricization, an SVD and a reshape, see \cite[Algorithm 1]{Oseledets2011}. Approximating a tensor for a given tolerance in the same sense as \eqref{eqn:eps_X} can also be done in the TT-decomposition, see \cite[Theorem 2.2]{Oseledets2011}.

The last requirement was that we need to be able to calculate a TT-decomposition of $\Ten{W}_Z$ (the equivalent of $\VM W_Z$ in \eqref{eqn:A(X)}), directly from the decomposition of $\Ten{Z}$. To see this, we recall the case of $d = 1$. Let $\VM 1_q$ be a vector with all ones of length $q$. 

%We first note that $$\VM Z \cdot \text{diag}([f_i(\omega_1), f_i(\omega_2), \hdots, f_i(\omega_{n_1})]) = \VM Z ._H ( \VM 1_{n_0} [f_i(\omega_1), f_i(\omega_2), \hdots, f_i(\omega_{n_1})])$$ from which it follows that 
%$$\VM W_Z = \left( (\VM 1_m \otimes \VM I_{n_0}) \VM F \right) ._H \left( (\VM I_{n_0} \otimes \VM 1_m) \VM Z \right).$$ This means that in case of multiple parameters, $\Ten{W}_Z$ can be generalised to
%$$ \Ten{W}_Z = \left( \VM I_m \otimes \VM 1_{n_0}, \VM I, \hdots, \VM I \right).\Ten{F} ._H \left( \VM 1_m \otimes  \VM I_{n_0}, \VM I, \hdots, \VM I \right).\Ten{Z}.$$
%From \cite{Oseledets2011}, we obtain that 
%\begin{equation} \label{eqn:TT_W_Z}
%\left\{\begin{matrix*}[c]
%\VM G^{W_Z}_j(i_j) & = & \VM G^F_j(i_j) \otimes \VM G^Z_j(i_j) & & j = 1,2, \hdots, d, i_j = 1,2, \hdots, n_j, \\ 
%\VM G^{W_Z}_0 & = & \VM G^F_0 \otimes \VM G^Z_0 &  & j = 0.
%\end{matrix*}\right.
%\end{equation}
%\textcolor{red}{
%---------------------------------------------------------------------------------------------------------------------------------------------------------------------
%Hetgeen hierboven staat klopt niet. Dit is juiste(re) versie:
We first note that $$\VM Z \cdot \text{diag}([f_i(\omega_1), f_i(\omega_2), \hdots, f_i(\omega_{n_1})]) = \VM Z ._H ( \VM 1_{n_U} [f_i(\omega_1), f_i(\omega_2), \hdots, f_i(\omega_{n_1})])$$ from which it follows that 
$$\VM W_Z = \left( (\VM 1_{n_U} \otimes \VM I_m) \VM F \right) ._H \left( (\VM I_{n_U} \otimes \VM 1_m) \VM Z \right).$$
Note that we use here the left-kronecker product $\otimes$. This means that in case of multiple parameters, $\Ten{W}_Z$ can be generalised to
$$ \Ten{W}_Z = \left( (\VM 1_{n_U} \otimes \VM I_m), \VM I, \hdots, \VM I \right).\Ten{F} ._H \left( \VM I_{n_U} \otimes \VM 1_m, \VM I, \hdots, \VM I \right).\Ten{Z}.$$
From \cite{Oseledets2011}, we obtain that 
\begin{equation} \label{eqn:TT_W_Z}
\left\{\begin{matrix*}[c]
\VM G^{W_Z}_j(i_j) & = & \VM G^Z_j(i_j) \otimes \VM G^F_j(i_j) & & j = 1,2, \hdots, d, i_j = 1,2, \hdots, n_j, \\ 
\VM G^{W_Z}_0 & = & \VM G^Z_0 \otimes \VM G^F_0 &  & j = 0.
\end{matrix*}\right.
\end{equation}
%---------------------------------------------------------------------------------------------------------------------------------------------------------------------}
Note that the TT-decomposition is not unique. Moreover, it is possible that the ranks of the decomposition in \eqref{eqn:TT_W_Z} of the tensor $\Ten{W}_Z$ are overestimated.

If all functions $f_i, i = 1, 2, \hdots, m$ are separable, then there is no need to first calculate a TT-decomposition of $\Ten{F}$. We can get a decomposition directly from the separability of the functions $f_i$. As the functions $f_{i}$ are separable there are functions $g_{i_j}^{i}(\omega_{i_j})$ such that 
$f_{i} \left( \omega_{i_1}, \omega_{i_2}, \hdots, \omega_{i_d} \right) = g_{i_1}^{i}(\omega_{i_1}) g_{i_2}^{i}(\omega_{i_2}) \hdots  g_{i_d}^{i}(\omega_{i_d})$. It can be seen that a possible TT-decomposition is of the form
\begin{equation}
\left\{\begin{matrix*}[c]
\VM G^F_0 & = & \VM I_m \\
\VM G^F_j(i_j) & = & \text{diag}( [g_1^j(\omega_{i_j}), g_2^j(\omega_{i_j}), \hdots, g_m^j(\omega_{i_j})]) & & j = 1, 2, \hdots, d-1, i_j = 1,2, \hdots, n_j, \\
\VM G^F_d(i_d) & = & [g_1^d(\omega_{i_d}), g_2^d(\omega_{i_d}), \hdots, g_m^d(\omega_{i_d})]^T & & i_d = 1,2, \hdots, n_d.
\end{matrix*}\right.
\end{equation}
In this case, all interior ranks are $m$. This means that the interior ranks of the tensor $\Ten W_Z$ are maximal $m r^Z_i, i=0, 1, \hdots, d$. Note that also here the ranks of the decomposition are possibly overestimated. The lower are the ranks of the tensor $\Ten Z$, the more efficient the other operations are. 
To have this, it is important how the eigenvectors are stored in the tensor $\Ten{X}^\subs{V}$. Eigenvectors are only unique up to a complex constant, but for the ranks of the tensor $\Ten X^{\subs V}$ it makes a difference which constant is used. The best we can do is to compare the eigenvector $\VM x(\VM \omega)$ for a point $\VM \omega$ with a certain vector $\VM a$, and store instead of $\VM x(\VM \omega)$ the vector $\dfrac{(\VM a^H \VM x(\VM \omega))}{\norm{ \VM a^H \VM x(\VM \omega)}} \VM x(\VM \omega)$. Remark that if $\VM a$ and the eigenvector have no complex elements then this boils down to multiplying with $1$ or $-1$. A possible choice for $\VM a$ is to take the eigenvector corresponding to a sample point.

If there are not many points in each dimension, it might be better not to use tensors but to reshape the tensor $\Ten{X}^\subs{V}$ as a matrix and to apply the PRA algorithm as described in Algorithm~ \ref{algo:resid_Arn_param}. The advantage of using a decomposition of the tensor is to reduce the storage and to reduce the number of operations, but if the number of parameter values is low in each dimension, then we cannot gain from the TT-representation as the cost of decomposing is too high in comparison with the possible gain in the remaining steps of the algorithm.

\subsection{Numerical Experiments} \label{sect:Tensors_NumExp}
In this section we extend Example \ref{ex:conv_diff} to multiple parameters and give the convergence results. The code can be downloaded following \href{http://twr.cs.kuleuven.be/research/software/delay-control/TensorKrylov}{this link}.

\begin{figure}
\setlength{\figW}{7cm} % determines the width of the images
 \begin{subfigure}[b]{0.5\textwidth}
  %\hspace{0.8cm}
  \input{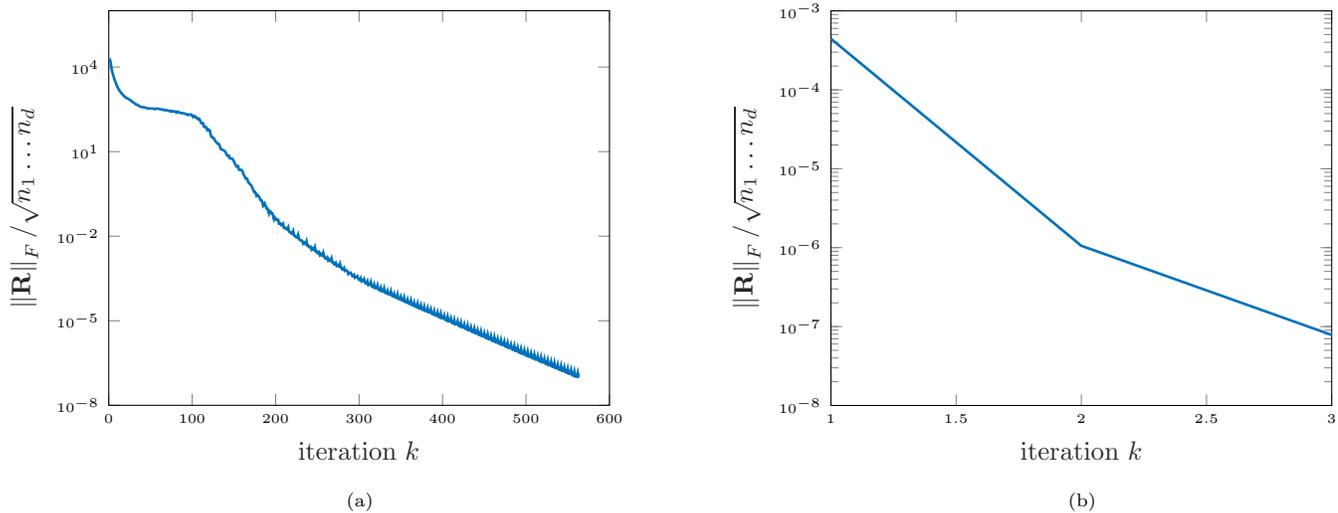}
  \caption{}
\end{subfigure} 
 \begin{subfigure}[b]{0.5\textwidth}
  %\hspace{0.8cm}
  % This file was created by matlab2tikz.
%
%The latest updates can be retrieved from
%  http://www.mathworks.com/matlabcentral/fileexchange/22022-matlab2tikz-matlab2tikz
%where you can also make suggestions and rate matlab2tikz.
%
%This file was created in matlab file C:\Users\Koen\Google Drive\Transcend\KULeuven\matlab_files\Low_rank Arnoldi\symmetrische matrix\fig_conv_diff_3D on 05-Feb-2020
%
\definecolor{mycolor1}{rgb}{0.00000,0.44700,0.74100}%
\begin{tikzpicture}

\begin{axis}[%
width=0.951\figW,
height=0.75\figW,
at={(0\figW,0\figW)},
scale only axis,
xmin=    1,
xmax=    3,
xlabel style={font=\color{white!15!black}},
xlabel={iteration $k$},
ymode=log,
ymin=1e-08,
ymax=0.001,
yminorticks=true,
ylabel style={font=\color{white!15!black}},
ylabel={$ \norm{\VM R}_F/\sqrt{n_1\hdots n_d}$},
axis background/.style={fill=white}
]
\addplot [color=mycolor1, forget plot, line width=1.0pt]
  table[row sep=crcr]{%
    1	0.000443512\\
    2	1.0556e-06\\
    3	7.83218e-08\\
};
\end{axis}

\begin{axis}[%
width=1.227\figW,
height=0.92\figW,
at={(-0.16\figW,-0.101\figW)},
scale only axis,
xmin=    0,
xmax=    1,
ymin=    0,
ymax=    1,
axis line style={draw=none},
ticks=none,
axis x line*=bottom,
axis y line*=left
]
\end{axis}
\end{tikzpicture}%
  \caption{}
\end{subfigure}
\caption{Results of the PRA algorithm with maximal subspace dimension $100$ for Example \ref{ex:conv_diff_3D}. In (a), we see the convergence when taking $n_1 = n_2 = n_3 = 10$. We see that we need quite a large number of iterations but tests showed us that  using Arnoldi for all parameter points separately needs approximately the same number of iterations. If we take the resulting subspace $\subs{V}$ and we use this as the initial subspace for $n_1 = n_2 = n_3 = 30$, we see that we just need a few iterations} \label{fig:conv_diff_3D} 
\end{figure}

\begin{eexample} \label{ex:conv_diff_3D}
Consider the convection-diffusion differential equation \eqref{eqn:conv_diff} in Example \ref{ex:conv_diff}. We consider the case where $d_{12} = 0$, we vary the variables $c_1$ and $c_2$ both in the interval $[4,6]$ and we let the variable $d := d_{11} = d_{22}$ vary in the interval $[-1.1, -0.9]$. In this way we can rewrite \eqref{eqn:conv_diff_A} into the form 
$$\VM A(c_1, c_2, d) = c_1 \VM A_1 + c_2 \VM A_2 + d \VM A_3.$$
The desired eigenvalue is the eigenvalue with smallest real part, in this case the same as the eigenvalue closest to zero. The matrices are of size $10^3 \times 10^3$ and we choose $n_1 = n_2 = n_3 = 30$ resulting in a total of $n_0 = 27 \times 10^{3}$ parameter values. As tolerances we choose $\eta_R = \dfrac{10^{-3}}{2+10^{-3}}$ and $\eta^{(k)}_X = 10^{-3} \norm{ \hat{\VM R}^{(k)}}_F/\sqrt{n_0}$
\end{eexample}
We first solve the problem discretised in each dimension by $10$ interior discretisation points. We see the result in Figure \ref{fig:conv_diff_3D}~(a). We observe that we need many iterations for convergence to the solution which is not surprising as using Arnoldi's method for all parameter points separately needs approximately the same number of iterations.  The obtained subspace $\subs{V}$ is then used as initial subspace when discretising each dimension in $30$. We see in Figure \ref{fig:conv_diff_3D} (b) that we only need a few iterations.

\section{Parametric shift-inverted Residual Arnoldi method}
in each iteration of the residual Arnoldi method, we applied the operator $\Ten A$ (see \eqref{eqn:diagAX}) or its tensor counterpart to Ritz vectors contaminated with error. In this way we are able to compute eigenvalues which are extreme in some sense. Users are often interested in the eigenvalue close to a certain shift $\sigma$. For this, the shift-inverted Arnoldi method can be used. For the ease of notation we assume that $\sigma = 0$. 

In this case we need to be able to efficiently compute linear systems with $\VM A(\VM \omega)$ for all parameter samples. The inverse of $\VM A(\VM \omega)$ can in general not be written in the form \eqref{eqn:VW_A}.  However if all matrices except one is of rank one, i.e $\VM A(\VM \omega) = f_0(\VM \omega) \VM A_0 + \sum_{i=1}^m f_i(\VM \omega) \VM A_i$ with $\VM A_i = \VM u_i \VM q_i^T, i = 1, \hdots, m$ with $m$ a small number, it can be proven via the Sherman-Morrison-Woodbury formula \cite{Hager1989, Woodbury1950} that for a vector $\VM q$, it holds that 

\begin{equation} \label{eqn:Ainv_subs}
    \VM A(\VM \omega)^{-1} \VM q \in \text{span}( \VM A_0^{-1} \VM q, \VM A_0^{-1} \VM u_1, \VM A_0^{-1} \VM u_2, \hdots, \VM A_0^{-1} \VM u_m).
\end{equation}
We see that there is only one basis vector of the subspace that depends on $\VM q$, while other vectors are the same for any $\VM q$ and for all parameter samples. Hence, one could keep all other vectors during all iterations in subspace $\subs{V}$.

It is possible to directly apply the PRA algorithm on $\VM A(\VM \omega)^{-1}$, however, in \cite{LeeStewart2007, Lee2007}, the authors proposed an adaptation of RA. In exact arithmetic, if we apply the RA algorithm on $\VM A^{-1}$, we add the residual vector $\VM r_{A^{-1}} = \VM A^{-1} \VM y - \dfrac{1}{\lambda} \VM y$ with $(\lambda, \VM y)$ a Ritz-pair of $\VM A$. The authors proposed to first calculate the residual $\VM r = \VM A \VM y - \lambda \VM y,$ to solve the system
\begin{equation} \label{eqn:system_A}
    \VM A \VM u = \VM r.
\end{equation} 
and to add the resulting $\VM u$ to the subspace. The authors in \cite{LeeStewart2007, Lee2007} proved that $\VM u$ is parallel to $\VM r_{A^{-1}}$. This has the advantage that solving \eqref{eqn:system_A} inexactly, produces an error proportional to $\VM r$. 

As the authors of \cite{LeeStewart2007, Lee2007} claim that also here the error proportional to the residual can be as high as $10^{-3}$ without harming the convergence, it is more appropriate to calculate the vector via \eqref{eqn:system_A}. Their method was called the Shift-and-invert Residual Arnoldi algorithm and similarly to the PRA algorithm, we extend this method to the parametric case. Instead of \eqref{eqn:system_A} we only need to calculate  
$$\VM A_0 \VM u(\VM \omega) = \VM r(\VM \omega)$$
upto a relative tolerance $\eta_R$, as we keep the other basis vectors in \eqref{eqn:Ainv_subs} in the subspace $\subs{V}$ in all iterations.
The tolerance for the relative error when solving \eqref{eqn:system_A} is denoted with $\eta_R$ and the tolerance for the low-rank approximation of the eigenvector over the parameter, is denoted with $\eta_X$. 

\subsection{Numerical experiments} \label{sec:PSIRA_numexp} 

The following examples arise from spectrum-based stabilization methods for time-delay systems \cite{Fenzi2017,Michiels2011}. Time-delay stability can be inferred from the spectral abscissa (real part of the rightmost eigenvalue) of an eigenvalue problem, discretizing the infinitesimal generator operator by collocation approach \cite{Breda2015}. Spectrum-based stabilization is based on minimizing the spectral abscissa as a function of controller parameter and favors situations where the spectral abscissa is not smooth, since minima of the spectral abscissa are often characterized by defective rightmost eigenvalues \cite{Fenzi2019}. In this context, our novel approach may lead to a significant contribution in the probabilistic spectrum-based stabilization \cite{Fenzi2017} since it permits to tackle larger scale problem with higher stochastic dimension. The code can be downloaded following \href{http://twr.cs.kuleuven.be/research/software/delay-control/TensorKrylov}{this link}.

\begin{figure}
%\end{subfigure} 
\setlength{\figW}{7cm} % determines the width of the images
 \begin{subfigure}[b]{0.5\textwidth}
  %\hspace{0.8cm}
  \input{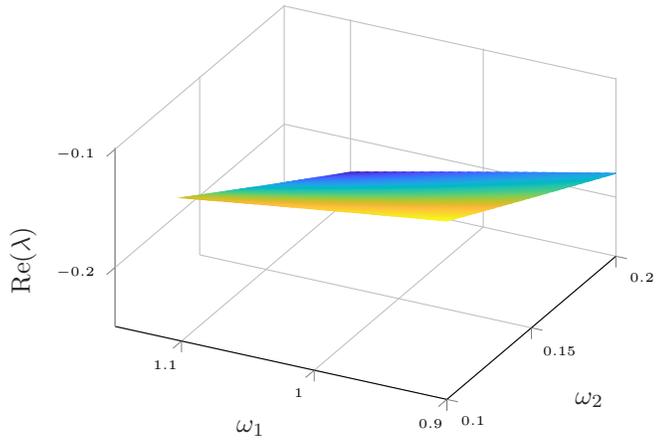}
  \caption{}
\end{subfigure} 
 \begin{subfigure}[b]{0.5\textwidth}
  %\hspace{0.8cm}
  \input{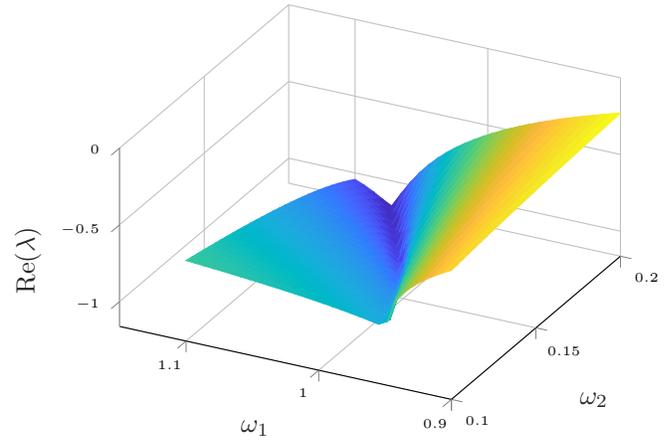}
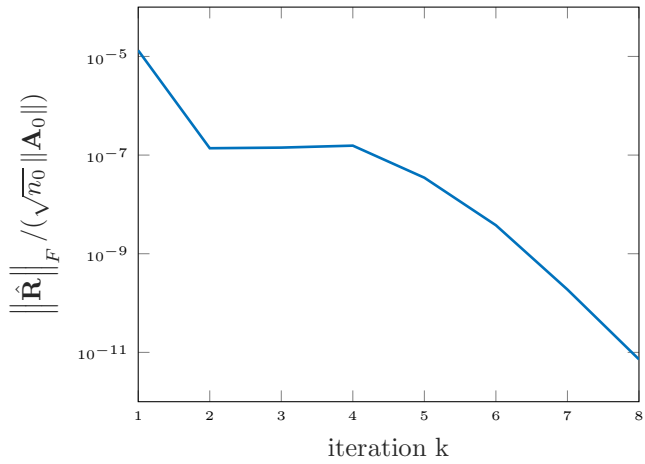
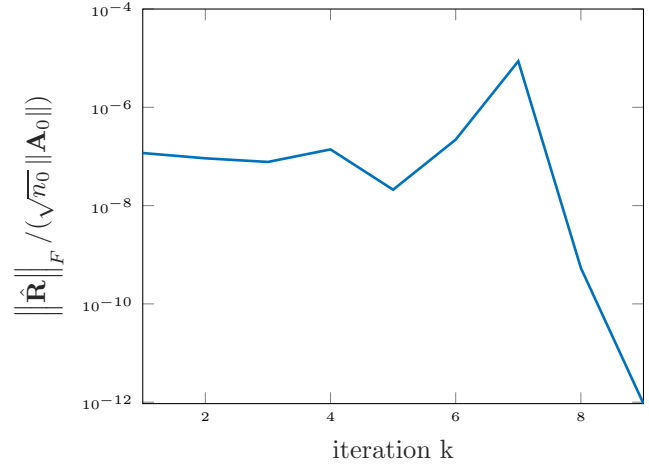
  \caption{}
\end{subfigure} \\
\setlength{\figW}{7cm} % determines the width of the images
 \begin{subfigure}[b]{0.5\textwidth}
  %\hspace{0.8cm}
  % This file was created by matlab2tikz.
%
%The latest updates can be retrieved from
%  http://www.mathworks.com/matlabcentral/fileexchange/22022-matlab2tikz-matlab2tikz
%where you can also make suggestions and rate matlab2tikz.
%
%This file was created in matlab file C:\Users\Koen\Google Drive\Transcend\KULeuven\matlab_files\Low_rank Arnoldi\Arnoldi_inverse\code_Luca\Delayed_Oscillator_Luca_main on 17-Feb-2020
%
\definecolor{mycolor1}{rgb}{0.00000,0.44700,0.74100}%
\begin{tikzpicture}

\begin{axis}[%
width=0.951\figW,
height=0.75\figW,
at={(0\figW,0\figW)},
scale only axis,
xmin=    1,
xmax=    8,
xlabel style={font=\color{white!15!black}},
xlabel={iteration k},
ymode=log,
ymin=1e-12,
ymax=0.0001,
yminorticks=true,
ylabel style={font=\color{white!15!black}},
ylabel={$\norm{\hat{\VM R}}_F/( \sqrt{n_0} \norm{ \VM A_0})$},
axis background/.style={fill=white}
]
\addplot [color=mycolor1, forget plot, line width=1.0pt]
  table[row sep=crcr]{%
    1	1.32295e-05\\
    2	1.37658e-07\\
    3	1.41762e-07\\
    4	1.55396e-07\\
    5	3.47643e-08\\
    6	3.78066e-09\\
    7	1.87833e-10\\
    8	7.23195e-12\\
};
\end{axis}
\end{tikzpicture}%
  \caption{}
\end{subfigure} 
 \begin{subfigure}[b]{0.5\textwidth}
  %\hspace{0.8cm}
  % This file was created by matlab2tikz.
%
%The latest updates can be retrieved from
%  http://www.mathworks.com/matlabcentral/fileexchange/22022-matlab2tikz-matlab2tikz
%where you can also make suggestions and rate matlab2tikz.
%
%This file was created in matlab file C:\Users\Koen\Google Drive\Transcend\KULeuven\matlab_files\Low_rank Arnoldi\Arnoldi_inverse\code_Luca\Delayed_Oscillator_Luca_main on 17-Feb-2020
%
\definecolor{mycolor1}{rgb}{0.00000,0.44700,0.74100}%
\begin{tikzpicture}

\begin{axis}[%
width=0.951\figW,
height=0.75\figW,
at={(0\figW,0\figW)},
scale only axis,
xmin=    1,
xmax=    9,
xlabel style={font=\color{white!15!black}},
xlabel={iteration k},
ymode=log,
ymin=9.34741e-13,
ymax=0.0001,
yminorticks=true,
ylabel style={font=\color{white!15!black}},
ylabel={$\norm{\hat{\VM R}}_F/( \sqrt{n_0} \norm{ \VM A_0})$},
axis background/.style={fill=white}
]
\addplot [color=mycolor1, forget plot, line width=1.0pt]
  table[row sep=crcr]{%
    1	1.17439e-07\\
    2	9.17174e-08\\
    3	7.73337e-08\\
    4	1.39158e-07\\
    5	2.10753e-08\\
    6	2.20424e-07\\
    7	8.68027e-06\\
    8	5.38298e-10\\
    9	9.34741e-13\\
};
\end{axis}
\end{tikzpicture}%
  \caption{}
\end{subfigure} 
\caption{The spectral abscissa over the parameters for the two different cases in Example \ref{ex:FbackCtrl1}. We see in (a) that the spectral abscissa varies smoothly over the parameters, in (b) this is not the case anymore. The convergence results for the two different cases in Example \ref{ex:FbackCtrl1} are plotted in (c) and (d). We see that in the two cases the algorithm converges within eleven iterations.} \label{fig:FbackCtrl1_eigw_conv}
\end{figure}

\begin{eexample} \label{ex:FbackCtrl1} 
The spectral abscissa function, associated to the oscillator with feedback delay system \cite[Example~2]{Fenzi2019},  is approximated by the Infinitesimal Generator approach  \cite{Breda2015} as the spectral abscissa of the following parametric matrix: 
\begin{equation}\label{eq:Aex2}
\VM A(\omega_1, \omega_2;K_1, K_2) = \VM A_0(K_1, K_2) -\omega_1^2  \VM A_1 +-2\omega_1\omega_2\VM A_2, \text{ with } \omega_1 \in [0.9, 1.1],\; \omega_2 \in [0.1, 0.2].
\end{equation}
The non-singular matrix $\VM A_0$ and the rank-$1$ matrices $\VM A_1$, $\VM A_2$ have dimension $1002 \times 1002$ and are given by:
\[
\VM A_0(K_1, K_2)=\begin{pmatrix}
0 & 1 & 0 & \cdots & 0 & 0 & 0\\
0 & 0 & 0 & \cdots & 0 & K_1 & K_2\\
  &   &   & \VM I_2 \otimes \VM D_M  &   &     &    
\end{pmatrix},\;
\VM A_1=\begin{pmatrix}
0 \\ 1 \\ 0 \\ \vdots \\ 0
\end{pmatrix}
\begin{pmatrix}
1 \\ 0 \\ 0 \\ \vdots \\ 0
\end{pmatrix}^T,\;
\VM A_2=\begin{pmatrix}
0 \\ 1 \\ 0 \\ \vdots \\ 0
\end{pmatrix}
\begin{pmatrix}
0 \\ 1 \\ 0 \\ \vdots \\ 0
\end{pmatrix}^T
\]
where $\VM D_M\in\mathbb{R}^{M\times M}$ is the differentiation matrix \cite[Chapter~6]{Trefethen2000} obtained with $M=500$ Chebyshev nodes as collocation points. The controller parameters $K_1$ and $K_2$ are selected so that the spectral abscissa function presents different behaviors on the parametric domain  $(\omega_1,\omega_2)\in[0.9, 1.1]\times [0.1, 0.2]$ as illustrated in Figures \ref{fig:FbackCtrl1_eigw_conv} (a) and (b): for $(K_1, K_2) = [0.2, 0.2]$  the spectral abscissa smoothly behaves, while for $(K_1, K_2) = [0.6179, -0.0072]$  the spectral abscissa presents non-Lipschitz continuous points, due to multiple non-semi-simple rightmost eigenvalues. We refer to \cite{Fenzi2019} for further information on the example. 

We use $\eta_R = \dfrac{ 10^{-3}}{2+10^{-3}}$ fixed and $\eta^{(k)}_X = 10^{-3} \norm{ \hat{\VM R}^{(k)}}_F/\sqrt{n_0}$ relative to the residual as explained in Section \ref{sect:Error_an}. The parameter intervals are both discretised into $30$ points, so we have a total of $n_0 = 900$ points for which the spectral abscissa needs to be found. The spectral abscissa is plotted in Figure \ref{fig:FbackCtrl1_eigw_conv} (a)-(b). The convergence results are shown in Figure \ref{fig:FbackCtrl1_eigw_conv} (c)-(d).
\end{eexample}

\begin{figure}
%\end{subfigure} 
\setlength{\figW}{7cm} % determines the width of the images
 \begin{subfigure}[b]{0.5\textwidth}
  %\hspace{0.8cm}
  \input{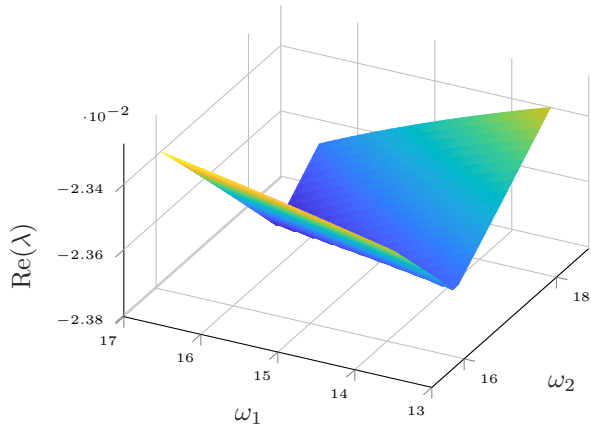}
  \caption{}
\end{subfigure} 
 \begin{subfigure}[b]{0.5\textwidth}
  %\hspace{0.8cm}
  \input{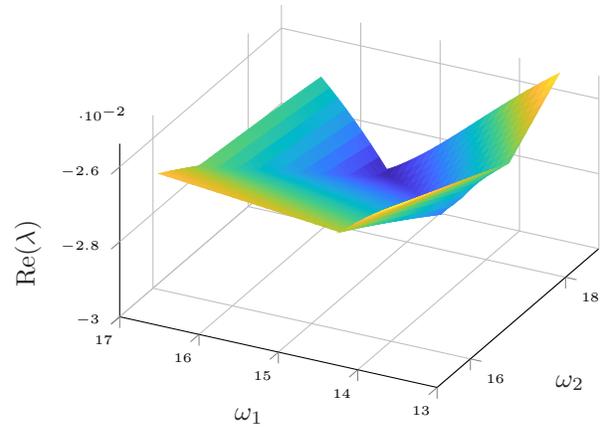}
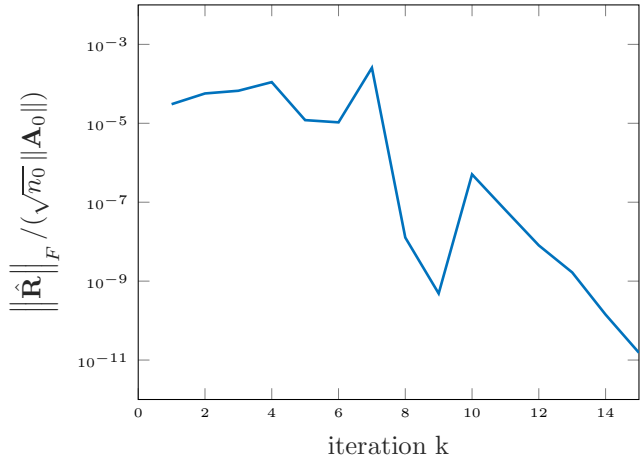
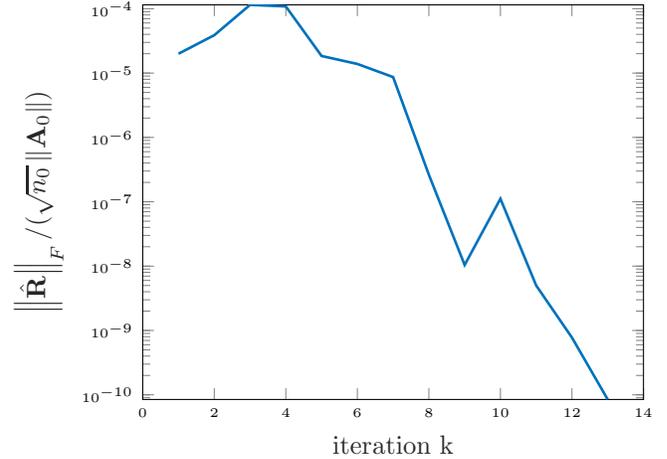
  \caption{}
\end{subfigure} \\
\setlength{\figW}{7cm} % determines the width of the images
 \begin{subfigure}[b]{0.5\textwidth}
  %\hspace{0.8cm}
  % This file was created by matlab2tikz.
%
%The latest updates can be retrieved from
%  http://www.mathworks.com/matlabcentral/fileexchange/22022-matlab2tikz-matlab2tikz
%where you can also make suggestions and rate matlab2tikz.
%
%This file was created in matlab file C:\Users\Koen\Google Drive\Transcend\KULeuven\matlab_files\Low_rank Arnoldi\Arnoldi_inverse\code_Luca\Heating_Luca_main on 17-Feb-2020
%
\definecolor{mycolor1}{rgb}{0.00000,0.44700,0.74100}%
\begin{tikzpicture}

\begin{axis}[%
width=0.951\figW,
height=0.75\figW,
at={(0\figW,0\figW)},
scale only axis,
xmin=    0,
xmax=   15,
xlabel style={font=\color{white!15!black}},
xlabel={iteration k},
ymode=log,
ymin=1e-12,
ymax= 0.01,
yminorticks=true,
ylabel style={font=\color{white!15!black}},
ylabel={$\norm{\hat{\VM R}}_F/( \sqrt{n_0} \norm{ \VM A_0})$},
axis background/.style={fill=white}
]
\addplot [color=mycolor1, forget plot, line width=1.0pt]
  table[row sep=crcr]{%
    1	3.04635e-05\\
    2	5.70015e-05\\
    3	6.67384e-05\\
    4	0.000110437\\
    5	1.20292e-05\\
    6	1.05661e-05\\
    7	0.000254667\\
    8	1.2734e-08\\
    9	4.86908e-10\\
   10	5.07167e-07\\
   11	6.31797e-08\\
   12	8.04593e-09\\
   13	1.66396e-09\\
   14	1.41777e-10\\
   15	1.53195e-11\\
};
\end{axis}
\end{tikzpicture}%
  \caption{}
\end{subfigure} 
\begin{subfigure}[b]{0.5\textwidth}
  %\hspace{0.8cm}
  % This file was created by matlab2tikz.
%
%The latest updates can be retrieved from
%  http://www.mathworks.com/matlabcentral/fileexchange/22022-matlab2tikz-matlab2tikz
%where you can also make suggestions and rate matlab2tikz.
%
%This file was created in matlab file C:\Users\Koen\Google Drive\Transcend\KULeuven\matlab_files\Low_rank Arnoldi\Arnoldi_inverse\code_Luca\Heating_Luca_main on 17-Feb-2020
%
\definecolor{mycolor1}{rgb}{0.00000,0.44700,0.74100}%
\begin{tikzpicture}

\begin{axis}[%
width=0.951\figW,
height=0.75\figW,
at={(0\figW,0\figW)},
scale only axis,
xmin=    0,
xmax=   14,
xlabel style={font=\color{white!15!black}},
xlabel={iteration k},
ymode=log,
ymin=8.46979e-11,
ymax=0.000114763,
yminorticks=true,
ylabel style={font=\color{white!15!black}},
ylabel={$\norm{\hat{\VM R}}_F/( \sqrt{n_0} \norm{ \VM A_0})$},
axis background/.style={fill=white}
]
\addplot [color=mycolor1, forget plot, line width=1.0pt]
  table[row sep=crcr]{%
    1	1.9977e-05\\
    2	3.87839e-05\\
    3	0.000114763\\
    4	0.000108584\\
    5	1.84359e-05\\
    6	1.38654e-05\\
    7	8.66993e-06\\
    8	2.62343e-07\\
    9	1.04244e-08\\
   10	1.11461e-07\\
   11	5.00292e-09\\
   12	7.76036e-10\\
   13	8.46979e-11\\
};
\end{axis}
\end{tikzpicture}%
  \caption{}
\end{subfigure} 
\caption{In (a)-(b), The spectral abscissa over the parameters for the two different cases in Example \ref{ex:Heating}. The convergence results for the two different cases in Example \ref{ex:Heating} can be seen in (c)-(d).} \label{fig:Heating_eigw_conv}
\end{figure}

\begin{eexample}  \label{ex:Heating}
% 	We consider the spectral abscissa function associated to an experimental heat transfer set up, consisting of 10 delay differential equations, which involves six different delays. The model has been described in \cite{Loiseau2009} and it is used as test cases for spectrum-based stabilization method in \cite{Michiels2010, Michiels2011, Fenzi2017}. 
	We analyze the spectral abscissa function associated to the heat transfer set up with two static feedback controllers $\VM K$, designed by probabilistic stabilization method considering the temperature on the left cooler and both cooler affected by uncertainty \cite[Section~5.3]{Fenzi2017}.  
	The spectral abscissa function is determined by the eigenvalues of the following $1030 \times 1030$ parametric matrix 
	\[\VM A(\omega_1, \omega_2; \VM K) = \VM A_0(\VM K) +  \dfrac{1-\omega_1}{\omega_1} \VM A_1 + \dfrac{1-\omega_2}{\omega_2} \VM A_2, \text{ where }\omega_1\in[13.5, 16.5],\, \omega_2\in[15.3, 18.7].\]
% 	The two static controllers $\VM K$, designed minimizing the mean of the spectral random variable modeling the uncertainty on the temperatures by independent uniform random variables  \cite[Section~5.3]{Fenzi2017},  lead to different behaviors of the spectral abscissa.  
	The two static controllers $\VM K$, designed in order to minimize the mean of the spectral abscissa, modelling the uncertainty on the temperatures by independent uniform random variables \cite[Section~5.3]{Fenzi2017},  lead to different behaviors of the spectral abscissa.  In Figure \ref{fig:Heating_eigw_conv} (a) resp. (b), the spectral abscissa for $ \VM K = [0.5105, -9.1810 \times 10^{-2}]$ resp. $\VM K=[0.45043,-0.1631]$ is used. The parameter intervals are both discretised into $30$ points, so we have a total of $n_0 = 900$ parameter values for which the spectral abscissa needs to be found. We use $\eta_R = \dfrac{ 10^{-3}}{2+10^{-3}}$ fixed and $\eta^{(k)}_X = 10^{-3} \norm{ \hat{\VM R}^{(k)}}_F/\sqrt{n_0}$ relative to the residual as explained in Section \ref{sect:Error_an}.  The convergence results are shown in Figure \ref{fig:Heating_eigw_conv} (c)-(d). 
\end{eexample}

We see in Figures \ref{fig:FbackCtrl1_eigw_conv} (c)-(d) and \ref{fig:Heating_eigw_conv} (c)-(d) that in both examples only a few iterations were needed. In Figure \ref{fig:FbackCtrl1_eigw_conv} (b) and \ref{fig:Heating_eigw_conv} (a)-(b) it can be observed that the spectral abscissa over the parameter is not differentiable as in some points (we only count eigenvalues with imaginary part equal or greater than zero) there are multiple eigenvalues for which the real part collides. We need to be careful when some eigenvalues of the large eigenvalues are close to each other. In our experience, subspace methods for calculating eigenvalues of non-normal matrices, can fail to converge, if the searched eigenvalue is close to other eigenvalues. Therefore, it is possible that due to the projection an eigenvalue is selected which is not the one we look for. The technique that we apply, is to save for every parameter sample in iteration $k$ five eigenvectors of the projected problem that best satisfies our criterium. In the folowing iteration $k + 1$ only an eigenvalue with associated eigenvector is selected if the angle between this eigenvector and one of the five selected eigenvectors in iteration $k$ is below a prescribed tolerance. 
 
\section{Conclusion}

In this paper, we proposed a method to calculate eigenvalues of a parametrized matrix $\VM A(\omega)$ for a number of parameter values at the same time. We started with introducing the Residual Arnoldi algorithm, which is a version of the Arnoldi algorithm that is robust against introduced errors. We applied the Residual Arnoldi Algorithm on a large number of sample points at the same time. To keep the computations feasible, we made low-rank approximations. On the one hand, we made a low-rank approximation of the searched eigenvector of the reduced eigenvalue problem over the parameters with a relative tolerance of $\eta_X$, and on the other, we made a low-rank approximation on the residual over the parameters with a relative tolerance of $\eta_R$. We saw that the total error relative to the residual needs to be smaller than $10^{-3}$, from which we deduce that $\eta_X$ needs to be chosen proportional to the residual and that $\eta_R$ can be chosen as high as $\dfrac{ 10^{-3}}{2+10^{-3}}$, see Section \ref{sect:conv_res_arn}.

We illustrated with numerical experiments that a higher $\eta_R$ even improves performance, since a higher tolerance leads to a slower increase of the subspace dimension and consecutely less restarts. At the end of the paper, we made the extension to mutltiple parameters and to shift-inverted Arnoldi for a specific range of problems. A possible extension is to make a version of the parametric shift-inverted Arnoldi method when more than one matrix has rank higher than one. The difficulty in this problem, is that we need to be able to solve parametrized linear systems efficiently. Another possibility for further research is the calculation of statistical moments of certain eigenvalues over the parameter.

\section*{Acknowledgements}
The authors thank the referees for their helpful remarks. We also thank Nick Vannieuwenhoven for some fruitful discussions about the used tensor format and we thank Luca Fenzi for the examples used in Section \ref{sec:PSIRA_numexp}. This work was supported by the project KU Leuven Research Council grant C14/17/072 and by the project G0A5317N of the Research Foundation-Flanders (FWO - Vlaanderen). This work does not have any conflicts of interest.

\bibliographystyle{elsarticle-harv}
\bibliography{ms}

\end{document}